\documentclass[journal ]{new-aiaa}
\usepackage[utf8]{inputenc}
\usepackage{textcomp}
\usepackage{subcaption,booktabs}
\usepackage{float}
\usepackage{graphicx}
\usepackage{amsmath}
\usepackage[version=4]{mhchem}
\usepackage{siunitx}
\usepackage{longtable,tabularx}
\title{Optimal Solutions to Deflect Earth Crossing Objects Using Laser}

\author{Vivek Verma \footnote{Doctorate student, Aerospace Engineering, Pennsylvania State University, USA; lrdvivekverma@hotmail.com.}, Shribharath B. \footnote{Shribharath B, Post doctoral fellow, Technion, Haifa, Israel; bshribharath92@gmail.com.}, and Mangal Kothari \footnote{Former Professor, Aerospace Engineering, Indian Institute of Technology Kanpur, India; mangalgnc@gmail.com}}

\begin{document}

\maketitle

\begin{abstract}
This paper solves and analyzes a trajectory optimization problem to deflect Earth-crossing objects (ECOs) employing continuous thrust obtained using a laser ablative system. The optimal control is determined for various initial ECO-Earth configurations to achieve the desired miss distance. The formulation incorporates the gravitational effect on the object due to the Earth using the patched-conic method. The constrained trajectory optimization problem is solved using Non-Linear Programming (NLP). First, the continuous control problem is solved, assuming both constant and variable power consumption, followed by a detailed comparison between the continuous control schemes. Subsequently, the work extends to studying sub-optimal solutions that can accommodate power fluctuations in the controller. The optimal control offers a range of alternative operational methods for asteroid deflection missions with trade-offs in power consumption and the total mission time. For impulsive deflection, the existing work reports two optimal solutions. One of the solutions is found to be better as it leads to a final ECO orbit that has its next Earth passage later than the other solution. Finally, the Moon's gravitational effect on the orbit of an ECO is studied. The reported results provide a comprehensive understanding of various scenarios in the process of ECO deflection.
\end{abstract}




\section{Introduction}\label{sec:intro}
\lettrine{N}{ear}-Earth asteroids, apart from being rich sources of minerals, are known to pose risks to life on the Earth through a catastrophic collision event. The continuous supply of asteroids from the main belt to Earth-crossing orbits aided by orbital resonances with Jupiter makes unsettling collision threats. Researchers have developed various deflection strategies to alter the orbit of an ECO to avoid collision with the Earth. The goal of these methods is to deflect the threatening object by changing its momentum by applying a small but timely external control force to achieve a safe separation between Earth and the asteroid during the close pass. The first step towards this is to estimate the change in orbital velocity required to avoid the collision. The work in 
\cite{doi:10.2514/3.11531} and 
\cite{ahrens_harris_1992} estimated the required change in orbital velocity by linearly approximating the orbital energy with respect to the velocity increment.  Further, Ref.~\cite{doi:10.2514/3.11531} proposed the usage of kinetic impact and nuclear explosion to provide the necessary change in velocity. Among the various possible deflection methods proposed,~\cite{park_mazanek_2003} has studied the effectiveness of different asteroid deflection methods and has concluded that ablation might be another efficient method. Ablation methods modify the trajectory of an object by irradiating it with laser light sources which results in continuous momentum chance. A laser with sufficient intensity can ablate the ECO's surface by causing a plasma blow-off. This process will generate a small thrust in the opposite direction to the applied laser changing the object's orbit. A deflection scheme based on the ablation method is presented in~\cite{McInnes1999} to generate the required thrust using solar radiation by attaching lightweight reflectors to the ECO in concern. Another technique utilizing solar power by concentrating the solar light in a small area using a large mirror to ablate the surface of the object has been proposed in~\cite{articleMaddock}. 
Moreover, last year NASA launched the DART mission to understand the behavior of the kinetic impact to deflect an asteroid. For this work, we assume a spacecraft carrying a laser ablation system along with necessary equipment is already deployed in the vicinity of ECO. 

A optimal control history is desired as it can lead to minimum control demand. A two-dimensional optimization problem to deflect ECOs using only two-body models was solved in~\cite{ParkandRoss1999}. They assumed the Sun-Earth and the Sun-asteroid systems to be independent neglecting the gravitational interaction between the asteroid and Earth. They found that the minimum impulse, $\Delta\text{v}$, required is a non-monotonous function of the impulse time (the time instant when impulsive control is applied) but has some local features (maxima and minima) over a secular variation that is inverse proportional to the impulse time (refer to Fig.~\ref{fig:thesis_Song_a12e06}). The finer structure is related to the orbital period of the object. The local minima occurred during the perihelion passes, making the earliest possible perihelion the ``optimal time'' for application of the $\Delta \text{v}$. A similar trend is also reported in~\cite{articleSong2010} between operation duration and operation start time for continuous control using laser ablation with constant power. In addition, Ref.~\cite{ParkandRoss1999} has pointed out that there exist two such optimal impulses that are equal in magnitude and opposite in direction (phase difference of $180^{\circ}$). This work was extended in~\cite{Ross2001} to study the Earth's gravitational effect on the asteroid during the close pass. The addition of Earth's gravity resulted in an increase in the required minimum $\Delta \text{v}$. Moreover, Ref.~\cite{Ross2001} observed that the increase in the magnitude of the impulse required was more pronounced for ECOs in a nearly circular heliocentric orbit with a size just about the same as Earth. The number of optimal solutions found remained two even after adding Earth's gravity. This work proposes a comparative study to identify a better solution between the two solutions. Based upon the analysis in~\cite{Ross2001}, the authors concluded that the problem of deflecting ECOs must include the gravity of Earth to avoid significant errors. The three-dimensional analysis of inclined ECO's orbits relative to the ecliptic plane done in~\cite{park_mazanek_2003} and~\cite{doi:10.2514/2.4814} showed that the optimal impulses required to deflection the ECOs are restricted to the orbital plane of the ECOs. With only a negligible component of $\Delta \mathbf{v}$ along the direction normal to the orbital plane. Generally, a simple analytical expression for the terminal condition of the deflection problem is achieved by analyzing it in the b-plane. In ~\cite{Ross2001},~\cite{park_mazanek_2003}, and~\cite{articleSong2010}, the modelling of the ECO's motion near the Earth have been done using three-dimensional patched conics. The b-plane is the plane formed by the velocity vector of asteroid relative to Earth and position vector of asteroid relative to Earth.

The known techniques to impart an impulse to the asteroid, such as \emph{kinetic impact} and \emph{nuclear explosion} may not precisely control the orbit. This is because these methods requires a single large impulse which when is imprecise may lead to the requirement of another large impulse to correct the errors caused. Further, the large impulses can lead to structural failure of the ECO producing too many fragments to be dealt with. The continuous control approach offers a more robust and precise technique to deflect an asteroid. Laser ablation is one such promising way to generate a  continuous thrust on the ECO. A pulsed laser ablative propulsion can also be utilized to alter the object’s trajectory by applying a series of small impulses. Generally, the pulse duration of the laser light is around in nanosecond or even less; hence the ablation process may be thought of as a continuous process. This technique is already investigated for the problem of deorbiting space debris using a ground-based laser facility in~\cite{1997ESASP.393..699B} and~\cite{articlephipps2003}. In~\cite{park_mazanek_2003}, linear analysis is done by assuming a linear relationship between the total impulse required and acceleration imparted by the laser ablation system to estimate the duration of operation of the laser ablation system based on the calculated minimum impulse. Later, assuming a continuous thrust from a laser source which has a constant power consumption, Ref.~\cite{articleSong2010} estimated the required operational time for the laser. It is found that the linear analysis done in earlier studies could yield up to $50\%$ error in the estimated total change in velocity required to attain the given miss distance. Moreover, the continuous control required with a variable power laser for control operation start time of less than one period of the ECO is studied in~\cite{articleSong2015}. For this study, it is assumed that a spacecraft carrying a laser ablation system has reached the vicinity of the asteroid and is capable of doing the required maneuvering. In the future, a spacecraft with VASIMR engines as the primary propulsion system is thought to be capable of performing rendezvous with an asteroid/comet belonging to near-Earth objects. VASIMR engines belong to the multi-megawatt class and are predicted to be available after the year 2050~\cite{ChangDiaz2000}. After successful rendezvous with the target ECO, the engine used to propel the spacecraft would also be used as a generator to supply the required power to the laser ablation system. Also, a study has been done in~\cite{ChangDiaz2000},~\cite{ParkYoung2005} to analyze a rapid rendezvous capability to reach target ECOs using future conceptual spacecraft with a VASIMR engine. Recently, in~\cite{GAMBI2021} the authors have discussed the swarm formations of spacecrafts around near-Earth asteroids to deflect them by ablating their surface using lasers.

The present work extends the study done in~\cite{park_mazanek_2003} and~\cite{articleSong2010} to find the control history for minimum power consumption utilizing a variable power laser. We also compared the results of constant power operation with variable power operation for laser operation times corresponding to more than one time period of the threatening asteroid. The above analysis is crucial for larger asteroids and robust deflection with a longer buffer time. In addition, the paper discusses several optimal and sub-optimal control solutions to minimize the total power consumption. The rest of the paper is organized as follows: Sec.~\ref{sec:2} discusses the goals and methods of the work in detail. Sec.~\ref{sec:3} presents the results of the numerical simulations and Sec.~\ref{sec:4} presents comparative studies for impulsive transfer methods. Finally, Sec.~\ref{sec:5} summarises the findings of this paper.

\section{Problem Formulation} \label{sec:2}
\begin{figure}[!ht]
	\centering
	\includegraphics[width=0.8\textwidth]{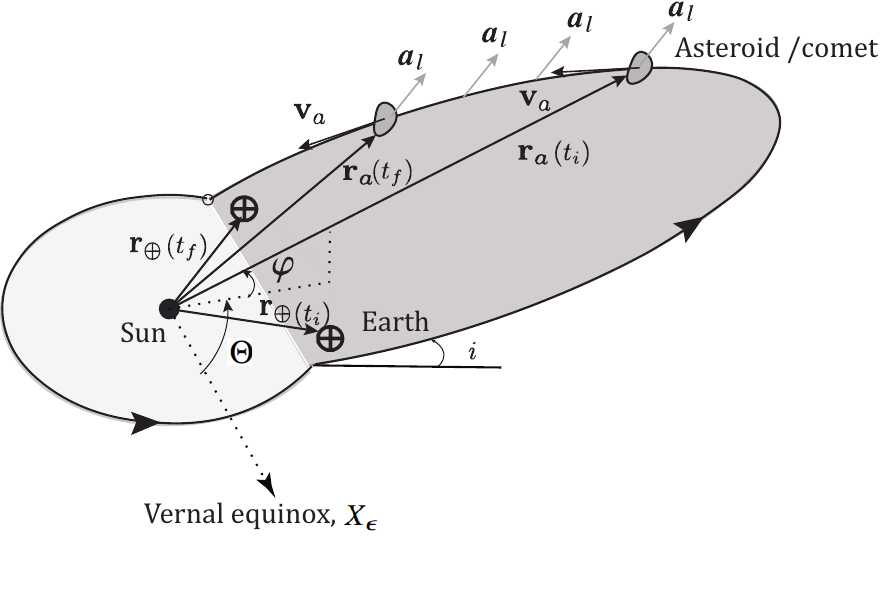}
	\caption{Schematic showing the orbit of bodies around the Sun (Not to scale)~\cite{articleSong2010}.}
	\label{orbitcontillus}
\end{figure}

The modeling and formulation of the problem are based on work in~\cite{articleSong2010} and~\cite{park_mazanek_2003}.
There are three bodies of interest: the Sun, Earth, and ECO. Other planets in the solar system are ignored for simplicity as the errors induced by them are usually negligibly small~\cite{park_mazanek_2003}. The Sun is taken to be stationary, while Earth and asteroid orbit the Sun, refer to Fig.~\ref{orbitcontillus}. In this work, the patched conic method is utilized to account for Earth's gravity on the asteroid.

\subsection{State dynamics}
The state equations describing the motion of ECO when ablation is active are formulated in the heliocentric spherical coordinates $\left\{\mathscr{S},\odot,\hat{\mathbf{e}}_r-\hat{\mathbf{e}}_{\Theta}-\hat{\mathbf{e}}_{\varphi}\right\}$; $\Theta$ is the angle between the $X_{\epsilon}-axis$ and the projection $\boldsymbol{r}$ in the ecliptic plane; $\varphi$ is the angle between the ecliptic plane and $\boldsymbol{r}$. In addition, a cylindrical coordinate system (CS) $\left\{ \mathscr{P},a,\hat{\boldsymbol{u}}-\hat{\boldsymbol{\zeta}}-\hat{\boldsymbol{w}} \right\}$ fixed  the center of the ECO in the perifocal plane; $\hat{\boldsymbol{u}}$ is a unit vector along the radial velocity direction, $\hat{\boldsymbol{\zeta}}$ is a unit vector along tangential velocity direction, which is perpendicular to $\hat{\boldsymbol{u}}$, and $\hat{\boldsymbol{w}}$ is a unit vector along the orbital angular momentum. Due to similarity, the unit vectors of spherical CS and cylindrical CS stay aligned at all times. The relation between the coefficients $(r_{X_{\epsilon}},r_{Y_{\epsilon}},r_{Z_{\epsilon}})^T$ of position vector $\boldsymbol{r}$ in heliocentric ecliptic CS and the coefficients $(r,\Theta,\varphi)^T$ in heliocentric spherical CS is
\begin{subequations}
	\begin{align}
		r_{X_{\epsilon}} &= r \cos \left(\varphi \right) \cos \left(\Theta \right) \\
		r_{Y_{\epsilon}} &= r \cos \left(\varphi \right) \sin \left(\Theta \right) \\
		r_{Z_{\epsilon}} &= r \sin \left(\varphi \right).
	\end{align}
\end{subequations}
Similarly, the relation between coefficients of the velocity vector $\mathbf{v}_a(t)$ \color{black} in the \color{black} cartesian CS and the spherical CS according to~\cite{articleSong2010} is
\begin{align}
	\left[\mathbf{v}_a(t)\right]_{uvw} = \begin{pmatrix}
		\cos \Theta \cos \varphi & \sin \Theta \cos \varphi & \sin \varphi\\
		-\sin \Theta & \cos \Theta  & 0\\
		-\cos \Theta \sin \varphi & -\sin \Theta \sin \varphi & \cos \varphi
	\end{pmatrix} \left[\mathbf{v}_a(t)\right]_{\epsilon}.
\end{align}
The set of governing equations that describes the motion of the ECO in the heliocentric spherical CS are~\cite{articleSong2010}
\begin{subequations}
	\begin{align}
		\dot{r} &= u\\
		\dot{u} &= \frac{v^2}{r} + \frac{w^2}{r} - \frac{\mu}{r^2} + a_l \sin \sigma \cos \beta\\
		\dot{v} &= -\frac{uv}{r} + \frac{vw \sin \varphi}{r \cos \varphi} + a_l \cos \sigma \cos \beta\\
		\dot{w} &= -\frac{uw}{r} + \frac{v^2 \sin \varphi}{r \cos \varphi} + a _l\sin \beta\\
		\dot{\Theta} &= \frac{v}{r \cos \varphi}\\
		\dot{\varphi} &= \frac{w}{r} ,
	\end{align}\label{decont}
\end{subequations}
where $r$ is the radial distance between the Sun and ECO. $u$, $v$ and $w$ are the components velocity $\mathbf{v}_a(t)$ in the radial, tangential, and the normal direction. The magnitude of acceleration imparted to the object by the surface jets caused due to laser ablation is denoted as $a_l$. The direction of the acceleration is specified using the angles $\sigma$ (in-plane angle) and $\beta$ (out-of-plane angle). $\sigma$ lies in the orbital plane and it is measured from $\hat{\textbf{v}}$ to the projection of $\boldsymbol{a}_l$ onto the perifocal plane. $(\beta)$ is the angle between the orbit plane and the thrust vector and is non-zero if the normal component of $\boldsymbol{a}_l$ is non-zero.

\subsection{The laser system}
A laser is a device that stimulates atoms or molecules to emit light at particular wavelengths~\cite{wikiLaser}. The laser light waves travel with their peak and trough all lined up in phase. In other words, the light from a laser device is coherent. The spatial coherence allows a laser to concentrate high energy on a small area. When a laser beam is emitted in pulses at a chosen frequency instead of a continuous wave, it is called a pulsed laser. The pulsed laser transfers the given energy with higher peak power than a continuous wave laser. Thus in a laser ablative propulsion, a pulsed laser is used to ablate the surface material of the ECO because it aids in quick heating of the material to a temperature that can cause sublimation. However, it is also possible to ablate material with a continuous wave laser sufficiently high in intensity.

The momentum coupling coefficient, $C_m$, is a common measure of the laser's energy consumption. It is defined as the ratio of impulse (momentum) generated to laser power. For a pulsed laser, it is the momentum generated per incident laser pulse energy.
\begin{equation}
	C_m = \frac{M \Delta \text{v}}{E_{laser}} = \frac{m v_e}{E_{laser}}, 
\end{equation}
where $M$ and $\Delta V$ are the mass and the velocity increment of the object ablated (ECO)and $m$, $v_e$ are the mass and the velocity of the exhausted propellant, and $E_{laser}$ is the laser energy delivered to the propellant. For propellants in solid-state, if the propulsion force is generated in vacuum, then, $C_m$ can be expressed as
\begin{equation}
	C_m = \frac{F_{av,laser}}{P_{av,laser}},
\end{equation}
where $F_{av,laser}$ and $P_{av,laser}$ are the time \color{black} averaged force on the ECO and the time averaged power of the laser respectively. Ablative laser propulsion is more powerful than pure photon propulsion. Even if total reflection is assumed, $C_m$ for the pure photon propulsion is $6.7\times10^{-9}\,N/W$ which is significantly lower compared to that for the ablative propulsion ($100\,N/MW$ to $10\,kN/MW$)~\cite{phipps2007laser}. The expression for the magnitude of acceleration imparted by a pulsed laser on an object it hits is derived in~\cite{articlephipps2003}. The force applied on the target by the ablation jets is
\begin{equation}
	F = P_{laser}C_m.
\end{equation}
Then, the acceleration is given as
\begin{equation}
	a =\frac{ P_{laser}C_m}{M}.
\end{equation}
The rate of mass loss due to mass ejection 
\begin{align}
	\begin{aligned}
		\dot{M} = \frac{P}{Q^*},
	\end{aligned}
\end{align}
where $Q^* = \frac{2\eta_{AB}}{C_m^2}$, 
$Q^*$ is called specific ablation energy, i.e. the energy required to ablate unit mass of the target. $\eta$ is the ablation efficiency defined to be the ratio of the exhaust kinetic energy imparted to mass element $dm$ to the incident laser energy. $C_m$, $Q$, $\eta$ are determined through experiments; their values depend upon the target material, intensity and duration of a laser pulse. A plasma blow-off using a sufficiently intense laser pulse will generate a small thrust that can change the ECO’s orbit given long enough time. Following the work of~\cite{articleSong2010} and~\cite{park_mazanek_2003}, we assume a laser system with power as high as $50\,MW$ in this work. Such futuristic energy generators have been used commonly in the literature of asteroid deflection.

\subsection{Problem Statement}
\begin{figure}[!ht]
	\centering
	\includegraphics[width= .7\textwidth]{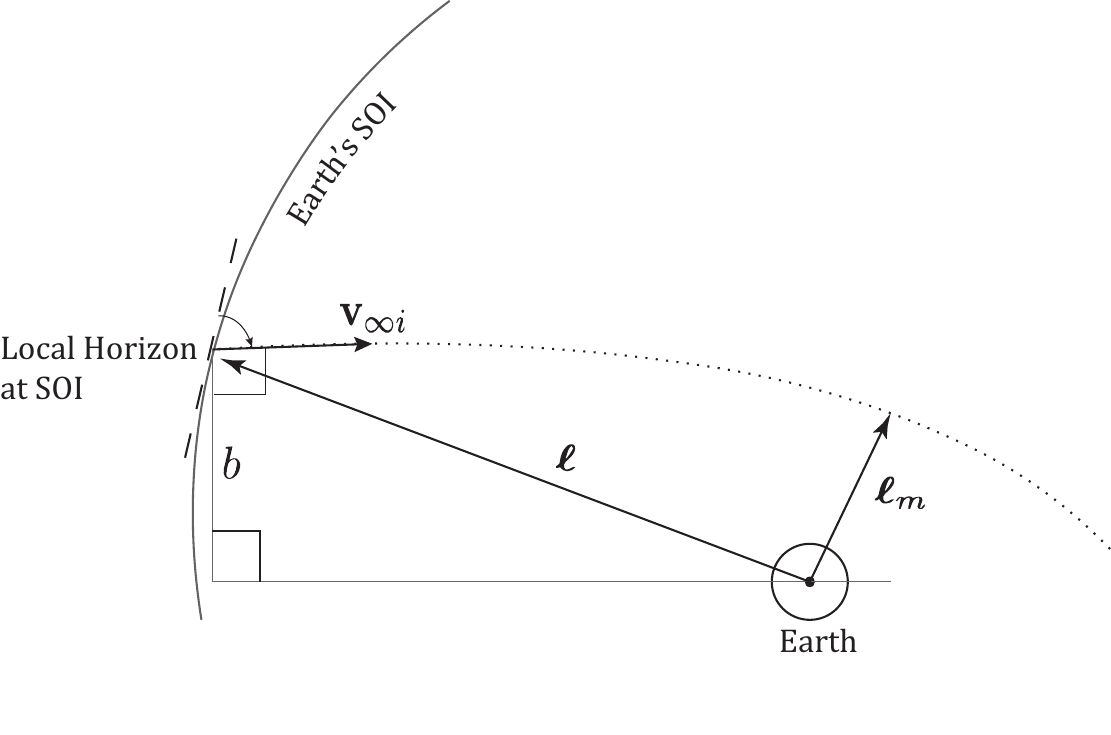}
	\vspace*{-1em}
	\caption{Schematic showing approach distance $b$ at the Earth's SOI in the $b-plane$ (Not to scale)~\cite{park_mazanek_2003}. $b-plane$ is defined by the vectors $\mathbf{v}_{\infty i}$ and $\boldsymbol{\ell}$. }
	\label{bplane}
\end{figure}

To begin with, a fictitious ECO is considered with orbital parameters such that it collides with Earth after some time. Then we set out to estimate the required continuous control history to deflect an ECO in its collision course with Earth using a pulsed laser ablation system. With the goal to deflect ECO with minimum control effort (which is equivalent to minimizing the energy consumed by the laser system) is subjected to the constraints. The objective function is defined as
\begin{equation}\label{objfuncont}
	J\left(\mathbf{u}(t),t_i, t_f\right) = \int_{t_i}^{t_f} \left|\mathbf{u}(t)\right|\,dt,
\end{equation}
where $t_i$ is the time instant when the laser ablation is started and $t_f$ is the time instant when it is stopped. The time interval $t_i-t_f$  is referred to as the operational time (abbreviate as $t_{op}$).
The terminal conditions are derived using the patched conic method. This method assumes that Earth's gravitational field affects the asteroid only inside the Earth's sphere of influence (SOI). Outside this region, the ECO orbit is affected solely by the Sun's gravity. The radius of Earth’s SOI is assumed to be $\ell_{soi}=9.31\times 10^8\ m$ as considered in~\cite{park_mazanek_2003}. Inside Earth's SOI, the ECO follows a hyperbolic trajectory as the relative velocity is higher than the escaping speed. Since to derive the terminal conditions, we have assumed a conic section (Hyperbola) orbit for ECO around the Earth. Hence, we refrain from applying the thrust to ECO inside the Earth's SOI. Moreover, the fraction of time ECO is inside the Earth's SOI is negligible relative to the total operation time. Figure~\ref{bplane} shows a fraction of the asteroid's orbit around Earth. The time instant at which the ECO enters the SOI region is $t_{soi}$ and for $t > t_{soi}$ the ECO will follow a hyperbolic orbit which has Earth at the nearest foci. The boundary condition of the optimal control problem is formulated in terms of the state of the ECO at $t= t_{soi}$. Let $\ell$ denotes the distance between Earth and the asteroid, $\dot{\ell}$ be the time derivative of $\ell$. The approach distance $b$ is the perpendicular distance between the Earth's center and the relative velocity vector of the asteroid relative to the Earth at the periphery of the Earth's SOI, refer to Fig.~\ref{bplane}. From the theory of conic sections, one can determine the relationship between the approach distance $b$ and the minimum separation between Earth and the ECO in the subsequent motion. At a critical value of the approach distance $b=b_i$, the resulting hyperbolic trajectory leads to exactly the same separation as desired. Hence, the terminal constraint to the optimal control problem is set that $b$ should be equal to the impact parameter $b_i$ at $t=t_{soi}$ when $\ell$ becomes equal to $\ell_{soi}$ so that we can achieve the desired miss distance $\ell_m$. Also, $\dot{\ell}$ is required to be negative, implying that the ECO should be approaching Earth at the boundary of the SOI. This condition helps the solver to accurately determine the value of $t_{soi}$. In addition to these, there are two more constraints. First, the $t_{op}$ should be positive, and the second is that the laser operation must end before the asteroid enters the Earth's SOI. The operation time is the time interval for which the laser ablative system is active and it is denoted as $t_{\text{oper}} = t_f - t_i$. Here, $t_i$ is the operation start time and $t_f$ is the time at which the operation stops.
In short, the terminal boundary conditions on relative distance and time are
\begin{align}\label{bdrycond}
	\begin{aligned}
		\ell_{soi}-\ell = 0,\\
		b-b_i=0, \\
		\dot{\ell}<0, \\
		t_{op} > 0,\\
		t_{soi} - t_{op} > 0.
	\end{aligned}
\end{align}

The relative distance $\ell$ is given by
\begin{align}
	\ell = \sqrt{\boldsymbol{r}_\oplus \cdot \boldsymbol{r}_\oplus + \boldsymbol{r}_a \cdot \boldsymbol{r}_a - 2(\boldsymbol{r}_\oplus \cdot \boldsymbol{r}_a)}.
\end{align}
Then, we have
\begin{align}
	\dot{\ell} = \left[\boldsymbol{r}_\oplus \cdot \mathbf{v}_\oplus + \boldsymbol{r}_a \cdot \mathbf{v}_a - \left(\boldsymbol{r}_\oplus \cdot \mathbf{v}_a + \boldsymbol{r}_a \cdot \mathbf{v}_\oplus  \right) \right]/\ell 
\end{align}
The subscript $\oplus$ identifies a quantity associated with Earth. The subscript $a$ is used for a quantity associated with the ECO, and the subscript $\odot$ shows the quantities associated with the Sun. $\boldsymbol{r}_{\oplus}$ is the position vector of Earth from the Sun's center, $\mathbf{v}_{\oplus}$ is the velocity vector of Earth with respect to the Sun, $\boldsymbol{r}_{a}$ is the position vector of the ECO from the Sun's center, and $\mathbf{v}_{a}$ is the velocity vector of the ECO with respect to the Sun. The approach distance, $b$, depends upon the radius of SOI and the orientation of the relative inbound velocity, $\mathbf{v}_{\infty i}$, of the ECO. The relative inbound velocity is the vector difference of the ECO's velocity and Earth's velocity at $\ell = \ell_{soi}$. The elevation angle, $\phi$ (see Fig.~\ref{bplane}), made by the velocity with the local horizon at SOI, satisfies the relation
\begin{equation*}
	\cos \left(\phi + \frac{\pi}{2}\right) = \frac{\mathbf{v}_{\infty i} \cdot \boldsymbol{\ell}_{soi}}{\text{v}_{\infty i}\ell_{soi}}.
\end{equation*}
The approach distance is given as
\begin{align}\label{brsoi}
	\begin{aligned}
		b = \ell_{soi} \cos \phi 
	\end{aligned} 
 \quad \text{or} \quad
    \begin{aligned}
	     b = \left|\left(-\boldsymbol{\ell}_{soi}\right) - \left(-\boldsymbol{\ell}_{soi} \cdot \frac{\mathbf{v}_{\infty i}}{\text{v}_{\infty i}}\right)\right|.
    \end{aligned}
\end{align}
The expression for impact parameter, $b_i$, in terms of the miss distance $\ell_m$ is given by~\cite{battin1987introduction}
\begin{equation}\label{birm}
	b_i = \ell_m \sqrt{1 + \frac{2\mu_{\oplus}}{\text{v}_{\infty i}^2 \ell_m}}.
\end{equation}


\section{Solution Approach}\label{sec:3}
The optimal trajectory problem is described by Eqns (\ref{objfuncont})\text{--}(\ref{decont}), and (\ref{bdrycond}). The next step is too transcribed these equations to form a nonlinear programming (NLP) problem. This transcription is done here using the shooting method~\cite{betts2010practicalch4}. First the time, $t$ is discretized into $N$ intervals with $N+1$ node points ($t_0, t_1,..., t_N$). Thus, the transcribed NLP problem is to

\textit{Minimize}
\begin{equation} \label{contobj}
	J(\mathbf{u}_0,...\mathbf{u}_N) =  \sum_{k = 0}^{N-1} \frac{h_k}{2}\left( \left|\mathbf{u}\right|_k + \left|\mathbf{u}\right|_{k+1} \right),
\end{equation}
\vspace*{-1em}
such that following boundary constraints are satisfied
\begin{align}
	\begin{aligned}
		\ell_{soi}-\ell = 0\\
		b-b_i=0 \\
		\dot{\ell}<0 \\
		t_{op} > 0\\
		t_{soi} - t_{op} > 0
	\end{aligned}
\end{align}
and limits on decision variables are
\begin{align} \label{limitaccel}
	\begin{aligned}
		\mathbf{u}_{low}\le\mathbf{u}_k\le \mathbf{u}_{upp},\quad k=0,1,...,N.
	\end{aligned}
\end{align}
The above NLP problem is then solved using the MATLAB\textsuperscript \textregistered \ optimization tool box. The decision variables of the NLP problem are $t_{soi},t_{op},\mathbf{u}_0,\mathbf{u}_1,....,\mathbf{u}_N$.\\
The following assumptions are made to simplify the problem:
\begin{itemize}
	\item[1] The fictitious ECO will hit the Earth's center if no action is taken.
	\item[2] The shape of the fictitious ECO is a sphere, and it is homogeneous throughout and it does not rotate.
	\item[3] A spacecraft with a powerful, futuristic laser ablation 
	system, has already reached the ECO.
\end{itemize}
A similar analysis was carried out by~\cite{articleSong2010} in which they have assumed that the laser operates at constant power. It means that magnitude of the force imparted is constant throughout the operation time. The objective function minimized in~\cite{articleSong2010} was $t_{op}$. In their case, since the power feed to the laser was constant, minimizing the the operational time is equivalent to optimizing the total energy consumed. In our work, the power is allowed to vary from zero to a set maximum value. Our solution methodology is validated with the results from~\cite{articleSong2010}. We consider an Apollo-class asteroid with elements and size as given in Table~\ref{Tab:a12e02}.
\begin{table}[h] 
	\centering
	\caption{Specification of the properties related to a asteroid}
	\begin{subtable}[h]{0.25\textwidth}
		\begin{tabular}{|l|l|l|}
		\hline
		$a$ & $e$ & $i$\\
		($au$) &  &($deg$) \\
		\hline 
		1.2 & 0.6 & 0\\
		\hline 
		\end{tabular}
		\caption{Orbit elements}
		\label{tab:elements_p3}
	\end{subtable}
	\hfill
	\begin{subtable}[h]{0.3\textwidth}
		\begin{tabular}{|l|l|l|}
			\hline
			Mass & density & Size \\
			($kg$) & ($kg/m^3$) & $m$\\
			\hline 
			1.6e9 & 3e3 & 100\\
			\hline
		\end{tabular}
		\caption{Asteroid size}
		\label{tab:size_p3}
	\end{subtable}
	\hfill
	\begin{subtable}[h]{0.3\textwidth}
		\begin{tabular}{|l|l|l|}
			\hline
			Power & $C_m$ & $\eta_{AB}$ \\
			($MW$) & ($Ns/J$) & \\
			\hline
			$10$ & $5e-5$ & $1$\\
			\hline
		\end{tabular}
		\caption{Laser system}
		\label{tab:Laser_p3}
	\end{subtable}
	\label{Tab:a12e02}
\end{table}
\subsection{Numerical setup}
The scaling of variables for computational purposes is as follows. Distances are scaled by the length unit $(LU)$, $1\,LU = 1\, au$ and time is scaled by the time unit $(TU)$, $1\,TU= P/ (2\pi)$ where $P$ is the period of the Earth's orbit around the Sun. Thus, all quantities for measuring speed get scaled by the speed unit $(SU)$, $1\,SU = 1\,LU/1\,TU\approx21\, km/s$. Recall that $t_i$ is the time instant corresponding to the control action's start point. Here, the impact time is set to zero so that $t_i$ may be interpreted as the time left to impact if no action is undertaken. Also, recall that ~\cite{articleSong2010} and~\cite{doi:10.2514/2.4814} have demonstrated that the energy consumed normal to the orbital plane is negligible relative to the in-plane deflection energy required. Thus, the out-off-plane angle $(\beta)$ is taken to be zero. Then the decision variables are in-plane angle $(\sigma)$ and magnitude of the acceleration ($a_l$). The \textit{operation angle}, $\delta$, is defined to be the angle between $\mathbf{a}_l$ and $\mathbf{v}_a$ to visualize the orientation of control relative to the velocity of asteroid. To understand the merits and demerits of the variable magnitude control, it is necessary to contrast the solution obtained with constant power. Therefore, we first present results obtained after assuming a constant power from the laser.

\subsection{Constant power deflection}
\begin{figure}[H]
	\subfloat[$t_i = 1\,Tp$]{
		\label{fig:thesis_conti_a12e06_sect_orbit}
		\includegraphics[width=0.5\textwidth]{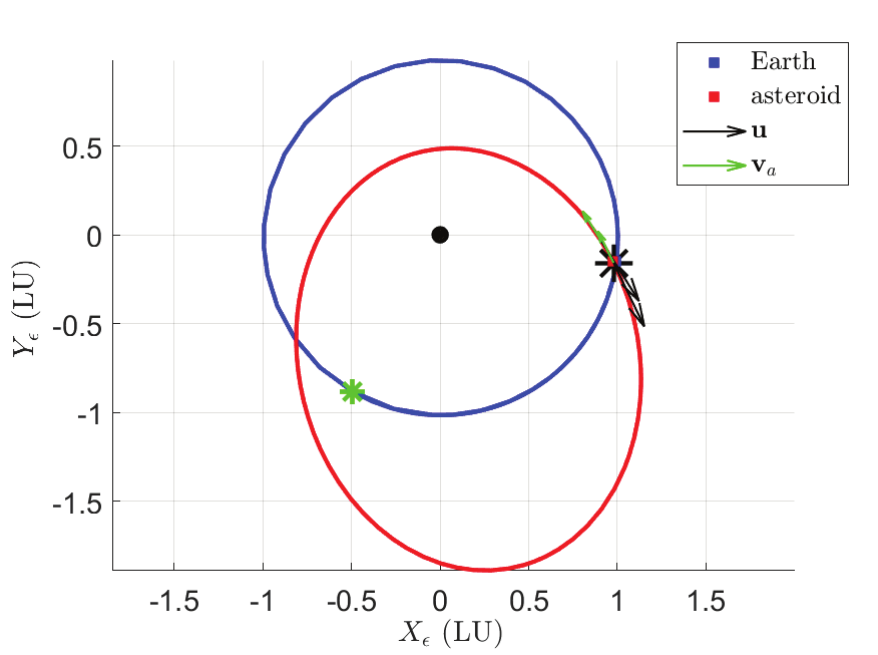} } 
	\subfloat[$t_i = 0.9\,Tp$]{
		\label{fig:thesis_conti_a12e06_angle}
		\includegraphics[width=0.5\textwidth]{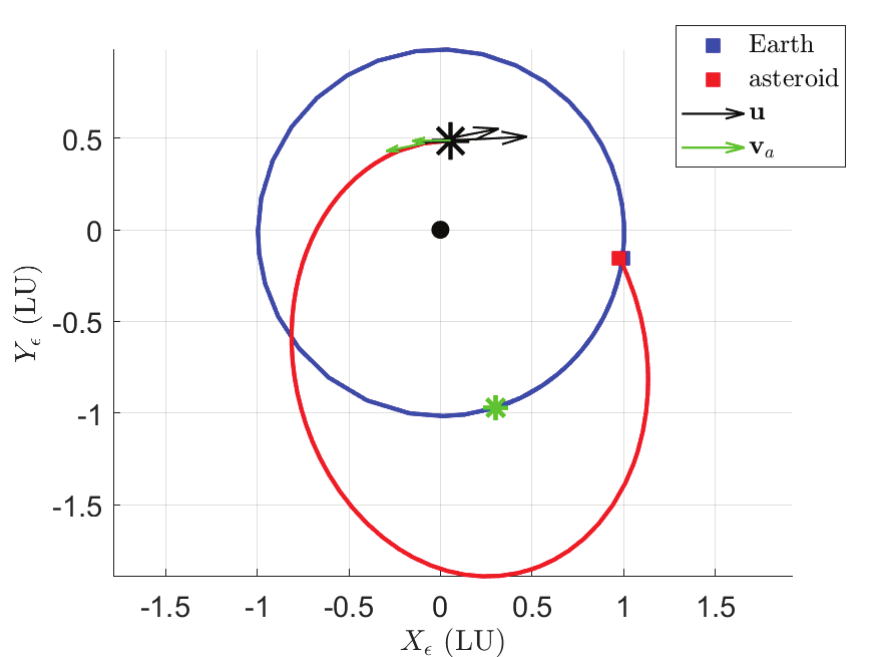} } 
	\caption{Showing orbit around the Sun. Both Earth and asteroid revolves in a counterclockwise direction. The black asterisk shows the starting point of the asteroid and the green asterisk shows the starting point of Earth.}
	\label{fig:thesis_conti_a12e06_orbit}
\end{figure}

\begin{figure}[h]
	\subfloat[$t_i = 1\,Tp$]{
		\label{fig:thesis_conti_a12e06_sect_orbit_090}
		\includegraphics[width=0.45\textwidth]{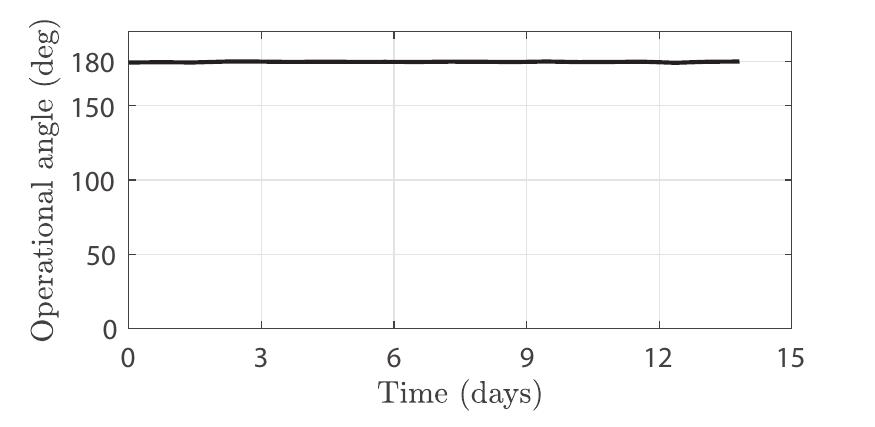} } 
	\subfloat[$t_i = 0.9\,Tp$]{
		\label{fig:thesis_conti_a12e06_angle_090}
		\includegraphics[width=0.45\textwidth]{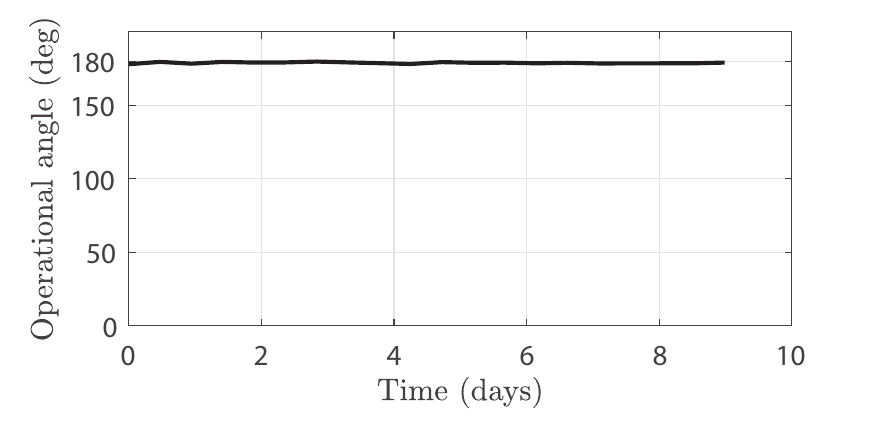} } 
	\caption{Illustrating the operational angle history (angle between the velocity and acceleration vector).}
	\label{fig:thesis_conti_angle_090}
\end{figure}

\begin{figure}[H]
	\centering
	\includegraphics[width= .6\textwidth]{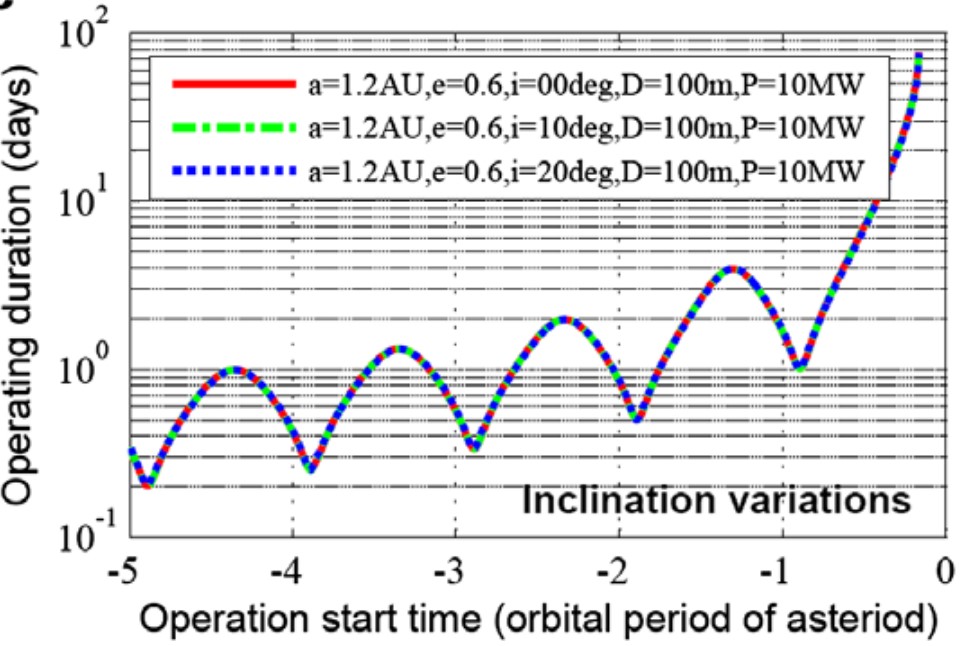}	  
	\caption{Minimum operation durations for various start times. Produced from~\cite{articleSong2010}.}
	\label{fig:thesis_Song_a12e06}
\end{figure}

The simulation results for constant power application are shown in Figs.~\ref{fig:thesis_conti_a12e06_orbit}\text{--}\ref{fig:thesis_Song_a12e06}. These figures illustrate the control history required for the start time $t_i = 1\,Tp$ and $t_i=0.9 \,Tp$, $1\,TP$ is equivalent to the time taken by the ECO understudy to complete one revolution (time period) around the Sun. Start time $t_i=0.9 \,Tp$ corresponds to the case when the asteroid is passing for the last time through perihelion before the collision. However, at $t_i=0.9 \, Tp$ asteroid is not precisely at perihelion but close to it. It is clear from Fig.~\ref{fig:thesis_conti_angle_090} that the optimal continuous acceleration required makes an angle of $180^0$ with the velocity vector of the asteroid in the orbital plane. This angle between velocity and required acceleration vector is termed as an \textit{operational angle}. For start time greater than the last perigee pass (for present scenario it is $0.90\,Tp$), the operational angle remains within the range $150^0$ to $210^0$~\cite{articleSong2010}.



\vspace*{-1em}
\begin{figure}[H]
	\subfloat[$t_i = 1\,Tp$\quad$t_{op}=13.82\,days$]{
		\label{fig:thesis_conti_a12e06_sect_orbit}
		\includegraphics[width=0.5\textwidth]{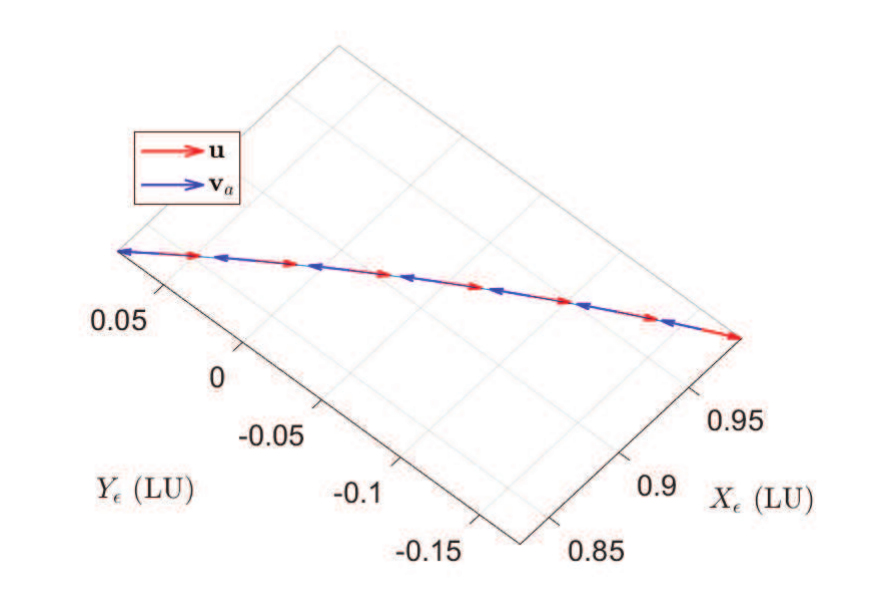} } 
	\subfloat[$t_i = 0.9\,Tp$\quad$t_{op}=8.98\,days$]{
		\label{fig:thesis_conti_a12e06_angle}
		\includegraphics[width=0.5\textwidth]{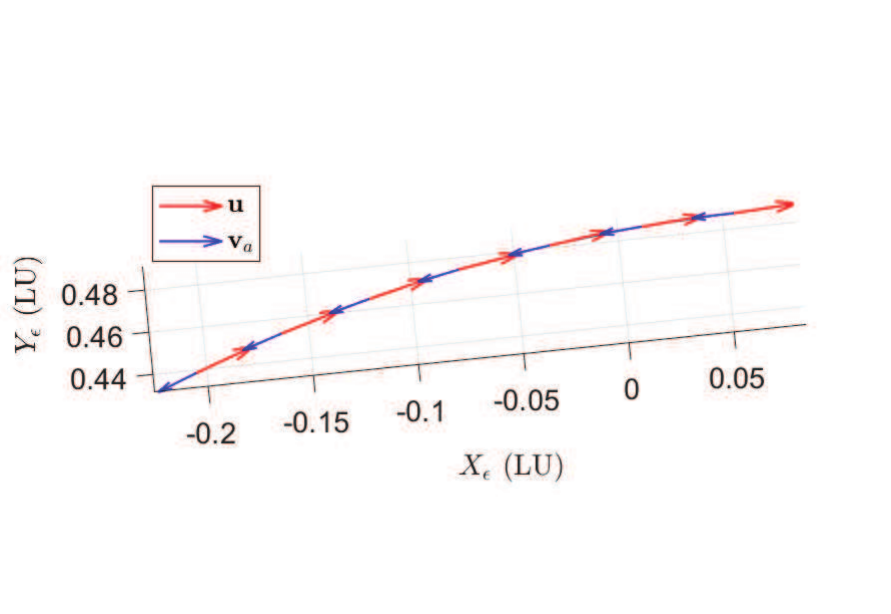} } 
	\caption{Illustrating the enlarge view of the continuous control in operation region.}
	\label{fig:thesis_conti_angle_1tp}
\end{figure}

The results for start times $t_i$ ranging from
$0.90$ to $1.9$ $T_p$ are listed in Table~\ref{Tab:timepower}. Figure~\ref{fig:thesis_Song_a12e06} shows the variation of the required operating time, $t_{op}$, with the operation start time, $t_i$, for the same operation setup as discussed in this work. Moreover, it illustrates that the curves corresponding to different inclination angles ($i=0^0,\ 10^0,\ 20^0$) coincide with each other. 
Here, we can see a locally cyclic behaviour of $t_{op}$ over a secular variation that increases with the start time. The series of local minima occur when the laser operation starts near to the perihelion. The \emph{oberth effect}, which states that the change in kinetic energy is proportional to the vehicle's velocity at the time of the burn can explain this feature of minima. The object's velocity reaches its peak when it passes through its periapsis. Therefore, the maximum change in kinetic energy occurs for a fixed change in velocity at the periapsis.

\subsection{Variable power} \label{sect_var_pow}
\begin{figure}[H]
	\subfloat[$t_i = 1.34\,Tp$]{
		\label{fig:thesis_conti_var_a12e06_sect_orbit}
		\includegraphics[width=0.5\textwidth]{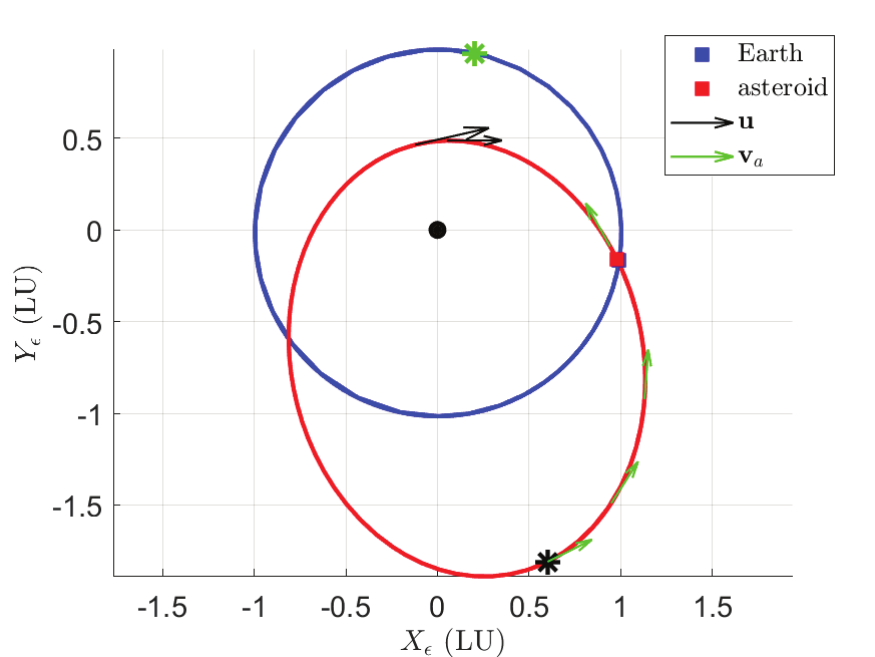} } 
	\subfloat[$t_i = 0.90\,Tp$]{
		\label{fig:thesis_conti_var_a12e06_angle}
		\includegraphics[width=0.5\textwidth]{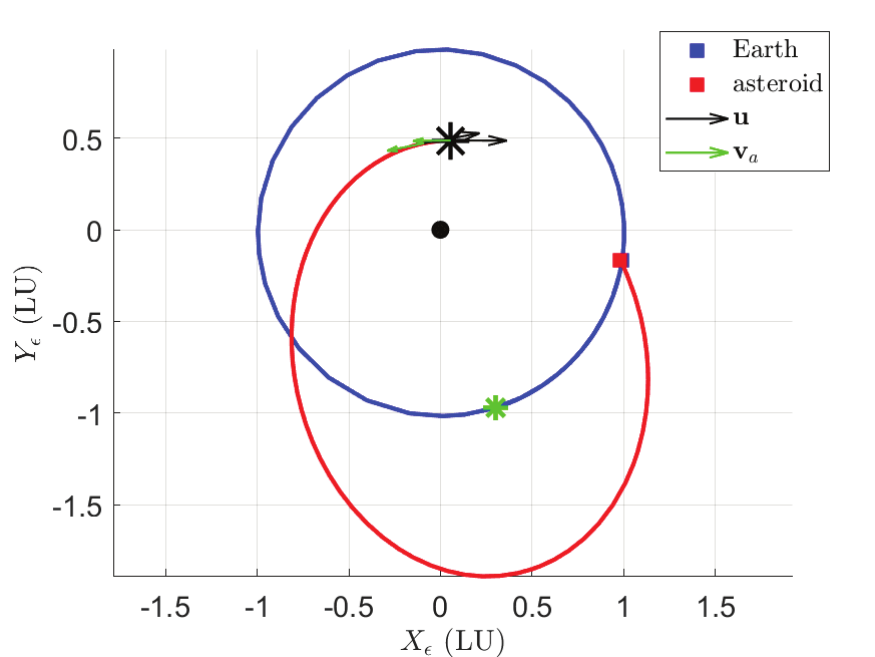} } 
	\caption{Showing orbit around the Sun. Both Earth and asteroid revolve in a counterclockwise direction. The black asterisk shows the starting point of the asteroid and the green asterisk shows the starting point of Earth.}
	\label{fig:thesis_conti_var_a12e06_orbit}
\end{figure}

\begin{figure}[H]
	\centering	
	\subfloat[$t_i = 1.34\,Tp$ ]{
		\label{fig:mag_1034}
		\includegraphics[width=0.5\textwidth]{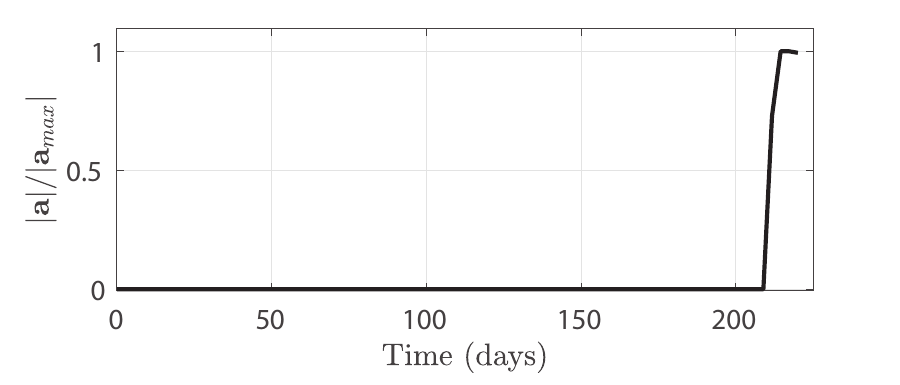} } 
	\subfloat[$t_i = 0.90\,Tp$]{
		\label{fig:angle_1034}
		\includegraphics[width=0.5\textwidth]{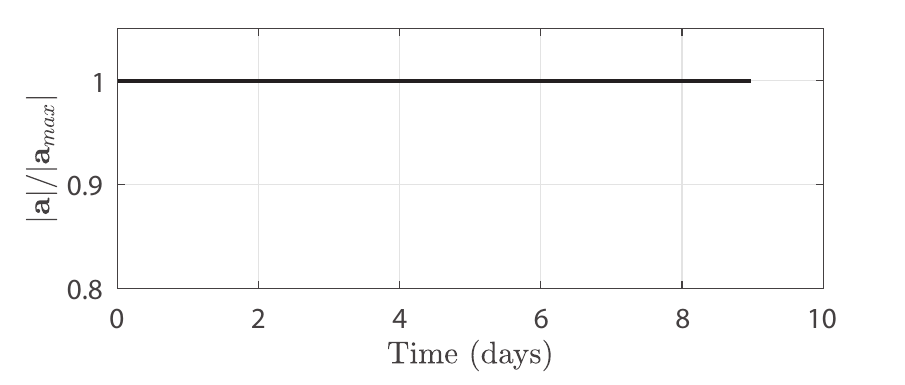} } 
	\caption{Illustrating the variation of acceleration magnitude of variable continuous control.}
	\label{fig:control_disc_1034}
	
\end{figure}

Consider the NLP problem defined by Eqns.~(\ref{contobj})--(\ref{limitaccel}) and the asteroid along with the laser configuration described in Table \ref{Tab:a12e02}. This time the decision variables are both in-plane angle $\sigma$ and magnitude of the acceleration $a_l$. Next, we will find the operation time and control history for various choices of operation start points and compare it with the results obtained using the constant power approach.
Figs.~\ref{fig:thesis_conti_var_a12e06_orbit}\text{--}\ref{fig:control_disc_1034} illustrate the results for variable power control operations. For start time $t_i = 0.90$, both variable and constant power operation yield the same results for optimal deflection. Another case with start time $t_i = 1.34$, when asteroid is near aphelion is presented to show the deviation in the behaviour of the variable power from the constant power operation. The control information for start time ranging $t_i=0.90$--$1.9$  
is listed in Table~\ref{Tab:timepower}.

\color{black}
\begin{figure}[H]
	\subfloat[$t_i = 1.34\,Tp$\quad$t_{op}=220.05$]{
		\label{fig:max_princ__vor}
		\includegraphics[width=0.5\textwidth]{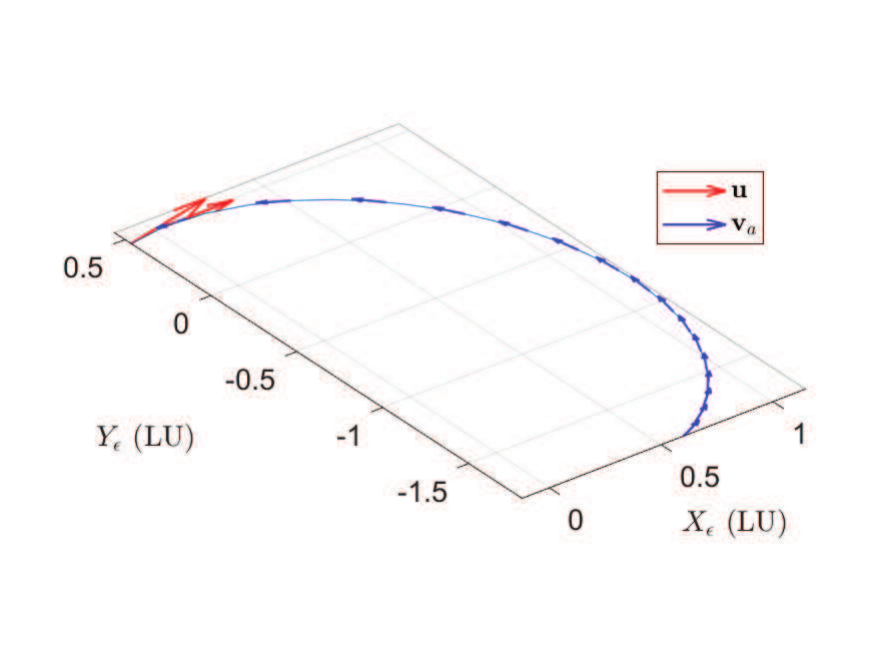} }
	\subfloat[$t_i = 0.90\,Tp$\quad$t_{op}=8.98$]{
		\label{fig:max_pring_nd}
		\includegraphics[width=0.5\textwidth]{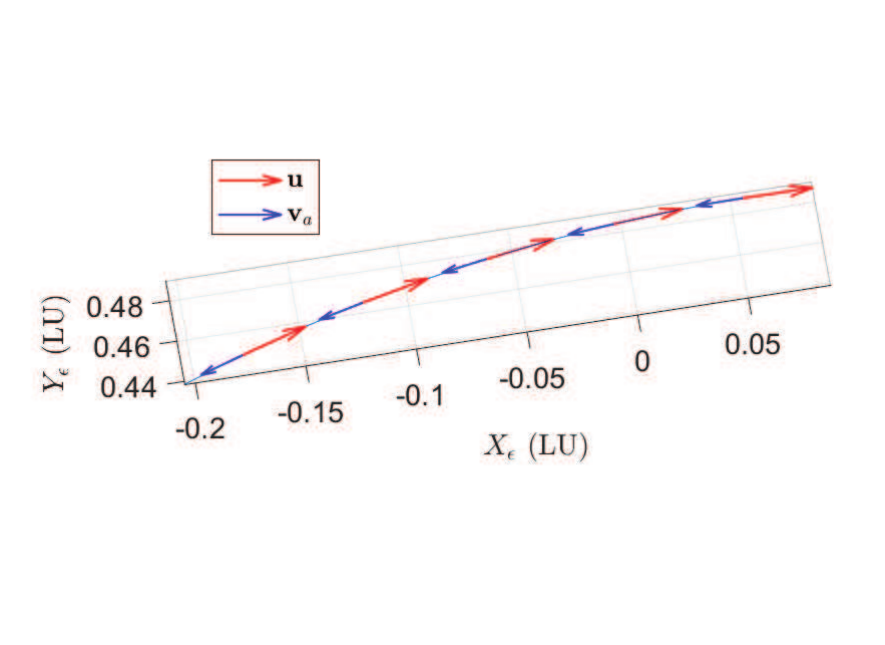} } 
	\caption{Illustrating the enlarge view of the continuous control in operation region.}
	\label{fig:comp_maxprin_pt2}
	
\end{figure}
\begin{figure}[H]
	\centering	
	\subfloat[$t_i = 1.34\,Tp$]{
		\label{fig:max_princ_10g_vor}
		\includegraphics[width=0.5\textwidth]{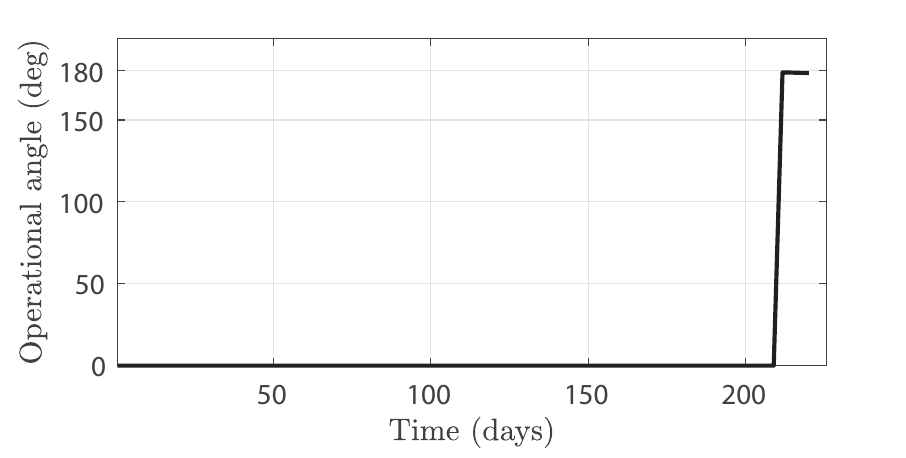} } 
	\subfloat[$t_i = 0.90\,Tp$]{
		\label{fig:max_prin_10g_nd}
		\includegraphics[width=0.5\textwidth]{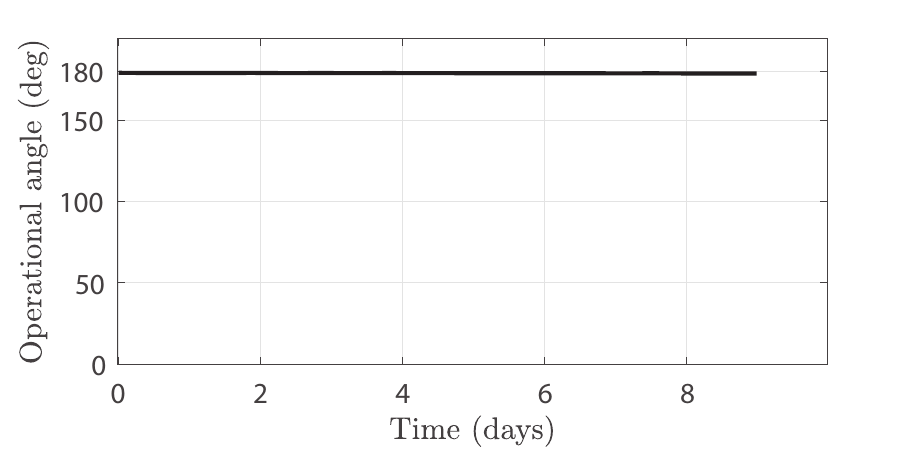} } 
	\caption{Illustrating the control angle for variable continuous operation.}
	\label{fig:thesis_angle_0895}
	
\end{figure}

\color{black}

It is seen from Table \ref{Tab:timepower} that the minimum energy consumption for the start time $t_i = 0.90\ Tp$ is almost the same for both the cases of variable and constant power consumption. However, as the start time increases, the total energy consumption remains almost the same until  $t_i= 1.70\ Tp$. While for constant power operation, the energy consumption increases up to $t_i= 1.34\ Tp$ and then decreases. Once again the minimum energy consumption for constant and variable power operation are similar for $t_i= 1.88\ Tp$ as the operation start time lies near the perihelion. In case of variable power operation, all the solutions ranging from $t_i= 0.90\ Tp$ -- $t_i= 1.70\ Tp$ are identical. Thus, for start times greater than $0.90\ Tp$ $(ti>0.90\ Tp)$, the laser operation does not start immediately and the laser system stays idle without operating for some time. We call this time gap as the \emph{idle time}. The laser is fired in the solution after the asteroid reached a region near \color{black} the perihelion and this gives an answer for the total energy consumed  similar to the one for $t_i = 0.90\ Tp$.

\begin{table}[h!]
	\centering                            
	\begin{tabular}{|l|l|l|l|l|l|}
		\hline
		  &\multicolumn{2}{|l|}{Constant power}&\multicolumn{3}{|l|}{Variable power}\\
		\hline
		Start time & Operation  &Energy&Idle&Actual operation & Energy \\
		Tp & time (day) & (kW day) &  time (day) &  time (day) & (kW day) \\
		\hline
		0.90$^1$& 8.98 & 89,846 & 0  & 8.99 & 89,836\\
		\hline 
		1& 13.82 & 1,38,153 & 51 & 9 & 89,876\\
		\hline 
		1.1& 20.54 & 2,05,380 & 96 & 10 & 89,945\\
		\hline 
		1.2& 26.63 & 2,66,337 & 156 & 10 & 89,860\\
		\hline
		1.34$^2$& 34.75& 3,47,455 & 209 & 11 & 89,835\\
		\hline
		1.5& 24.28  & 2,42,791  & 290 & 12 & 89,844 \\
		\hline
		1.7& 11.63  & 1,16,273  &  383& 12 & 89,806 \\
		\hline
		1.88$^3$& 4.49  & 44,855  &  0& 4.51 &44,658 \\
		\hline
		1.9& 4.50  & 44,979  & 3.7 & 5.2 & 44,542 \\
		\hline
	\end{tabular}
	\caption{Table enlist --- operation time and energy consumed for constant operation along with waiting time, operation time, and energy consumed for variable operation for various start times. 		1. local minima between $0$ - $1\,Tp$,
		2. local maxima between $1$ - $2\,Tp$,
		3. local minima between $1$ - $2\,Tp$, refer to Fig.~\ref{fig:thesis_Song_a12e06}.
		}
	\label{Tab:timepower}
\end{table}  

Additionally, Table~\ref{Tab:impulseNconti} compares the total change in velocity required by different control strategies namely, impulsive control and continuous control. For continuous control, both constant and variable power requirement are shown in the Table~\ref{Tab:impulseNconti}.

\begin{table}[h!]
	\centering                            
	\begin{tabular}{|l|l|l|l|}
		\hline
		  &Impulsive input& \multicolumn{2}{l|}{continuous input}\\
		\hline
		Start time &  & Constant & Variable \\
		\hline
		 & $\Delta |\mathbf{V}|$ & $\int_{t_i}^{t_f}|\boldsymbol{u}|\,dt$ &  $\int_{t_i}^{t_f}|\boldsymbol{u}|\,dt$ \\
		$Tp$ &  $m/s$ & $m/s$ &  $m/s$  \\
		\hline
		0.90& 25.01 & 24.71 &24.71 \\
		\hline 
		1& 41.39 & 37.99 & 24.72 \\
		\hline 
		1.1& 62.04 & 56.48 & 24.74\\
		\hline 
		1.2& 84.10 & 73.25 & 24.71 \\
		\hline
		1.34& 94.12& 95.56 & 24.71\\
		\hline
		1.5& 62.93  & 66.77  & 24.71 \\
		\hline
		1.7& 30.51  & 31.98  &  24.70 \\
		\hline
		1.88& 12.31  & 12.33  &  12.28  \\
		\hline
		1.9& 12.34  & 12.37  & 12.25 \\
		\hline
	\end{tabular}
		\caption{Table enlist --- the total change in velocity imparted by different control strategies.}
	\label{Tab:impulseNconti}
\end{table}

\subsection{Sub-optimal solutions}
\begin{figure}[H]
	\centering	
	\subfloat[Operational angle history ]{
		\label{fig:angle_history_1034_ll03}
		\includegraphics[width=0.5\textwidth]{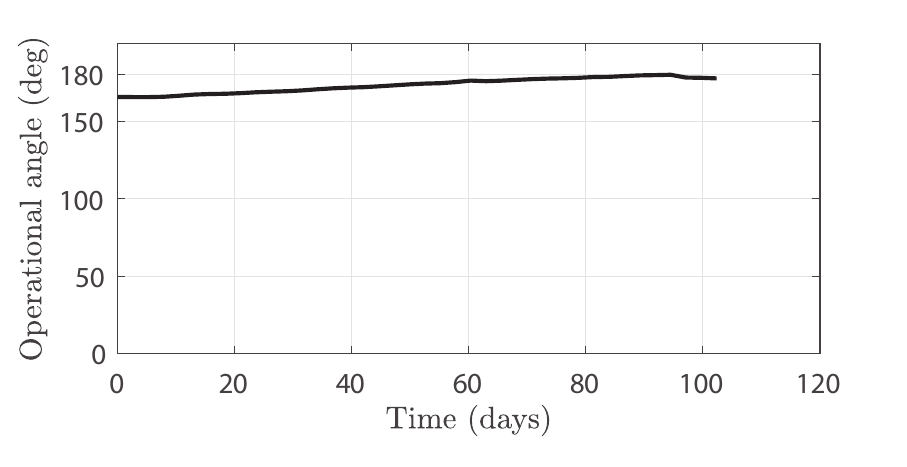} }
	\subfloat[Variation of acceleration magnitude]{
		\label{fig:mag_history_1034_ll03}
		\includegraphics[width=0.5\textwidth]{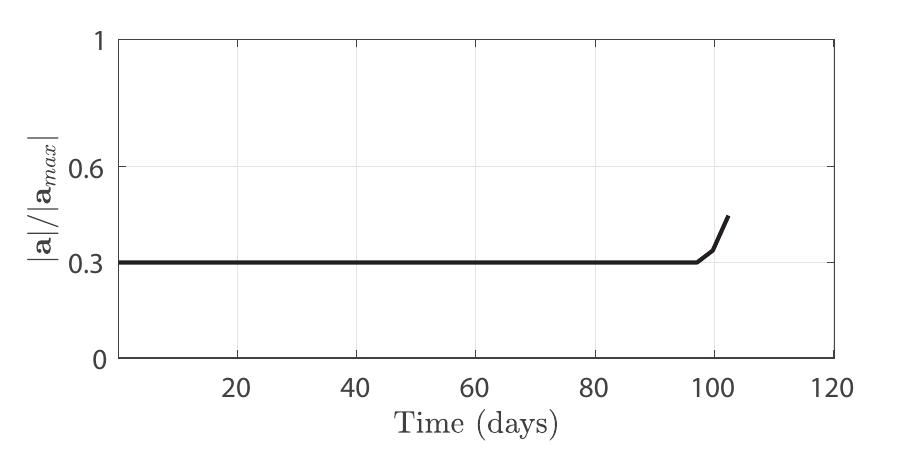} } 
	    \caption{Illustrating the continuous control required corresponding to start time $t_i = 1.34\,Tp$. The lower limit on the magnitude of the acceleration is $0.3$ times the maximum value at each node.}
	    \label{fig:thesis_subopti_control_1034_ll03}
\end{figure}

\begin{figure}[H]
	\centering	
	\subfloat[Operational angle history]{
		\label{fig:angle_history_1034_llinc}
		\includegraphics[width=0.5\textwidth]{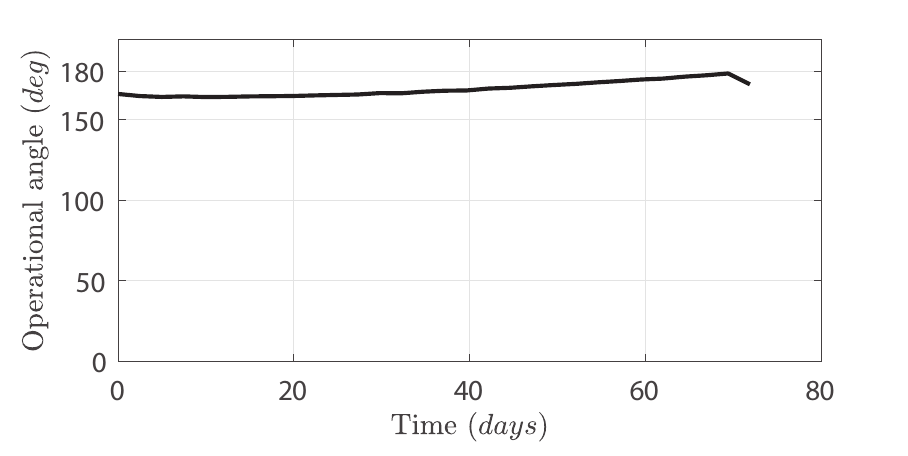} } 
	\subfloat[Variation of acceleration magnitude]{
		\label{fig:mag_history_1034_llinc}
		\includegraphics[width=0.5\textwidth]{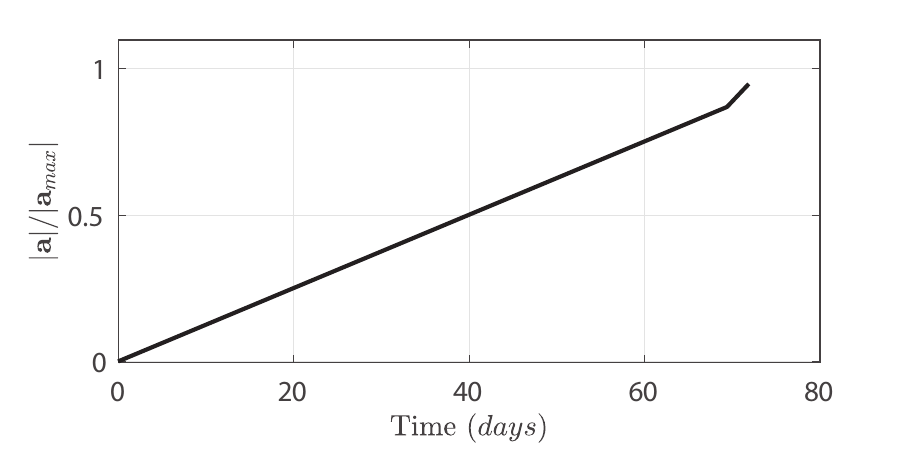} } 
	\caption{Illustrating the continuous control required corresponding to start time $t_i = 1.34\,Tp$. The lower limit on the magnitude of the acceleration is assumed to vary linearly from $0\,\text{--}\,0.9$ times the maximum value with time.}
	\label{fig:thesis_subopti_control_1034_llinc}
\end{figure}

\begin{figure}[H]
	\centering	
	\subfloat[Operational angle history]{
		\label{fig:angle_history_1034_lldec}
		\includegraphics[width=0.5\textwidth]{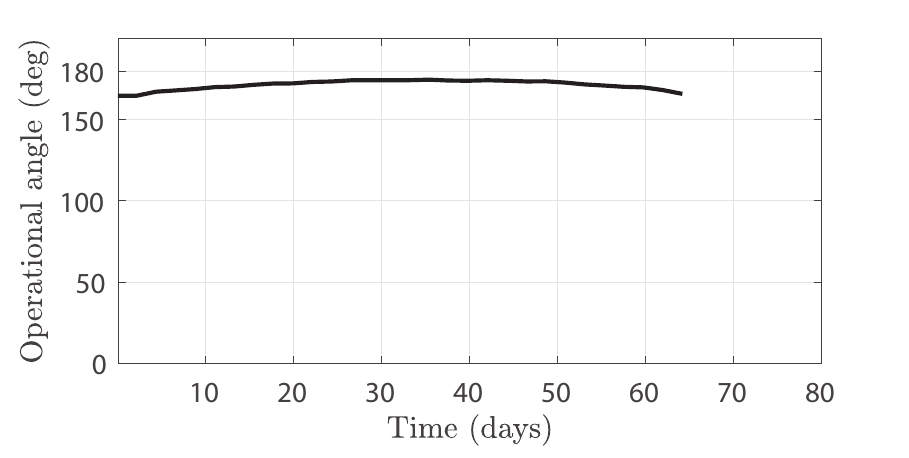} } 
	\subfloat[Variation of acceleration magnitude]{
		\label{fig:mag_history_1034_lldec}
		\includegraphics[width=0.5\textwidth]{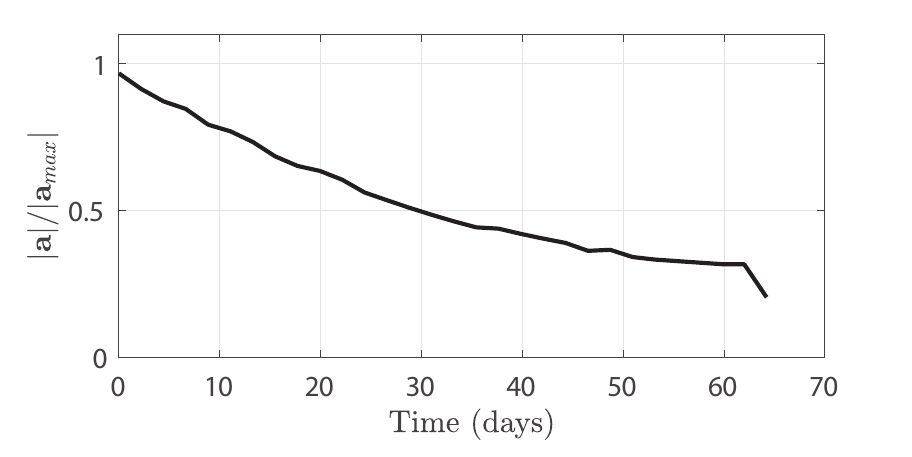} } 
	\caption{Illustrating the continuous control required corresponding to start time $t_i = 1.34\,Tp$.The lower limit on the magnitude of the acceleration is linearly varied from $0.9\,\text{---}\,0$ times the maximum value. }
	\label{fig:thesis_subopti_control_1034_lldec}
	
\end{figure}

This section will discuss sub-optimal solutions as we have seen in Section~\ref{sect_var_pow} that for the start time away from the perihelion has some idle time associated with it. The idle time depends upon how far the operation start point is from the coming perihelion point. And for this scenario, the variable power approach is unable to utilize the additional time in hand. One way to use the additional time available is to modify the control constraints, i.e., set a non-zero lower limit for the control magnitude. Figs.~\ref{fig:thesis_subopti_a12e06_sect_orbit_1034_ll03} -- \ref{fig:thesis_subopti_control_1034_lldec} show control history for variable operation for three different cases of lower bound on control magnitude. In case A, a constant lower bound equal to $0.3$ times the maximum possible acceleration of the laser ablation system is set, the results are depicted in Figs.~\ref{fig:thesis_subopti_a12e06_sect_orbit_1034_ll03} -- \ref{fig:thesis_subopti_control_1034_ll03}. Then, for case B, the lower bound is linearly varied from zero at $t_i$ to $0.9$ times maximum acceleration at $t_f$, shown in Figs.~\ref{fig:thesis_subopti_a12e06_sect_orbit_1034_llinc}--\ref{fig:thesis_subopti_control_1034_llinc}. Finally, for case C, the lower bound is set to vary linearly $0.9$ times maximum acceleration at $t_i$ to zero at $t_f$, illustrated in Figs.~\ref{fig:thesis_subopti_a12e06_sect_orbit_1034_lldec}--\ref{fig:thesis_subopti_control_1034_lldec}. Tables~\ref{Tab:thesis_subopti_t1}, \ref{Tab:thesis_subopti_t2}, and \ref{Tab:thesis_subopti_t3} present the operation intervals and energy consumption for the three cases respectively. From the data enlisted in Tables~\ref{Tab:timepower}, \ref{Tab:thesis_subopti_t1}--\ref{Tab:thesis_subopti_t3}, we can see that for the cases A, B, and C the energy requirement is less than the constant power operation. Along with this, the required operation time is less than the variable power with no limit on lower bound (i.e. lower limit is zero). It must be noted that this lower bound is purely computational and may not be a constraint of the real system. There is no change in the upper bound, and hence, the laser system used is the same as in earlier cases. 
In practice, it may be challenging to operate a system at its maximum capacity continuously and achieve the desired outcome. It may also be huge risk to wait until the last orbit of the ECO to deflect it as done in~\cite{articleSong2015}. These sub-optimal solutions can be beneficial in some cases, especially when the laser system may loses its full range of operation over the course of the deployment. One such scenario could occur when the mirror used by laser to concentrate the beam gets covered by ablative particles leading to a reduction in efficiency~\cite{GIBBINGS2013}.

\subsubsection{\textbf{Case} A}
\vspace*{-1em}
\vspace*{-2em}
\begin{table}[H]
	\begin{minipage}{0.4\linewidth}
		\caption{Specifying the results}
		\label{Tab:thesis_subopti_t1}
		\centering
		\begin{tabular}{|l|l|l|}
		\hline
		Start time & Operation  &Energy\\
		Tp & time (day) & (kW day)  \\
		\hline
		1.34& 102.34 & 3,09,935\\
		\hline 
	\end{tabular}
	\end{minipage}\hfill
	\begin{minipage}{0.6\linewidth}
		\centering
	    \includegraphics[width= 1\textwidth] {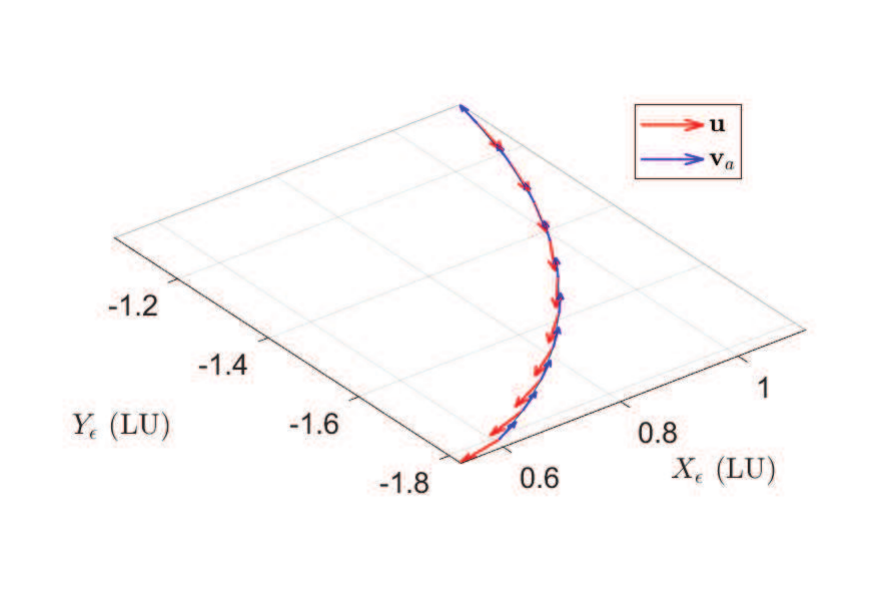}
	    \captionof{figure}{Enlarge view of the operation region.}
	    \label{fig:thesis_subopti_a12e06_sect_orbit_1034_ll03}
	\end{minipage}
\end{table}

\subsubsection{\textbf{Case} B}

\vspace*{-1em}
\begin{table}[H]
	\begin{minipage}{0.4\linewidth}
	\caption{Specifying the results}
		\label{Tab:thesis_subopti_t2}
		\centering
		\begin{tabular}{|l|l|l|}
		\hline
		Start time & Operation  &Energy\\
		Tp & time (day) & (kW day)  \\
		\hline
		1.34& 71.88 & 3,24,414\\
		\hline 
	\end{tabular}
	\end{minipage}\hfill
	\begin{minipage}{0.6\linewidth}
	   \centering
	   \includegraphics[width= 1\textwidth]     {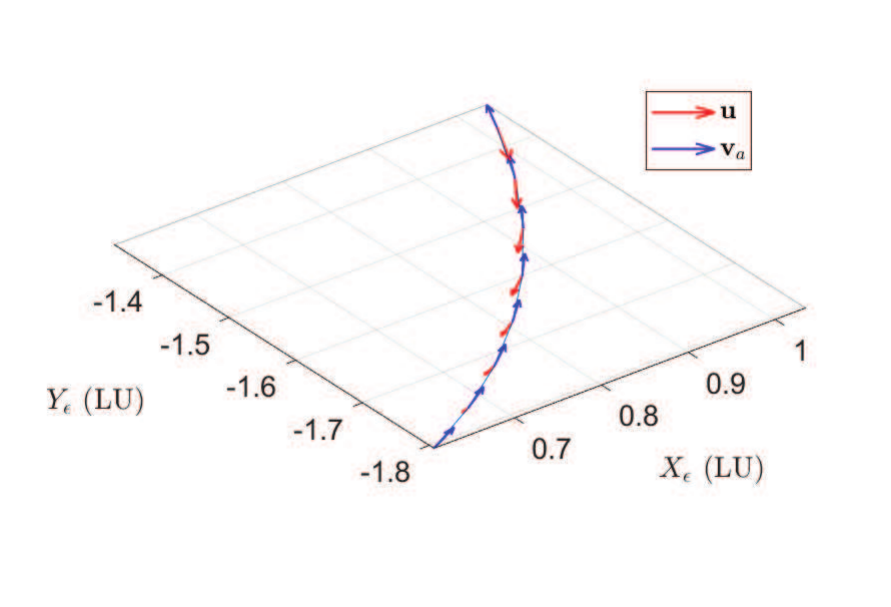}	  
	   \captionof{figure}{Enlarge view of the operation region.}
        \label{fig:thesis_subopti_a12e06_sect_orbit_1034_llinc}
	\end{minipage}
\end{table}
\vspace*{-2em}
\subsubsection{\textbf{Case} C}
\vspace*{-2em}
\begin{table}[H]
	\begin{minipage}{0.4\linewidth}
		\caption{Specifying the results}
		\label{Tab:thesis_subopti_t3}
		\centering
		\begin{tabular}{|l|l|l|}
		\hline
		Start time & Operation  &Energy\\
		Tp & time (day) & (kW day)  \\
		\hline
		1.34& 64.22 & 3,41,711\\
		\hline 
	\end{tabular}
	\end{minipage}\hfill
	\begin{minipage}{0.6\linewidth}
		\centering
	\includegraphics[width= 1\textwidth] {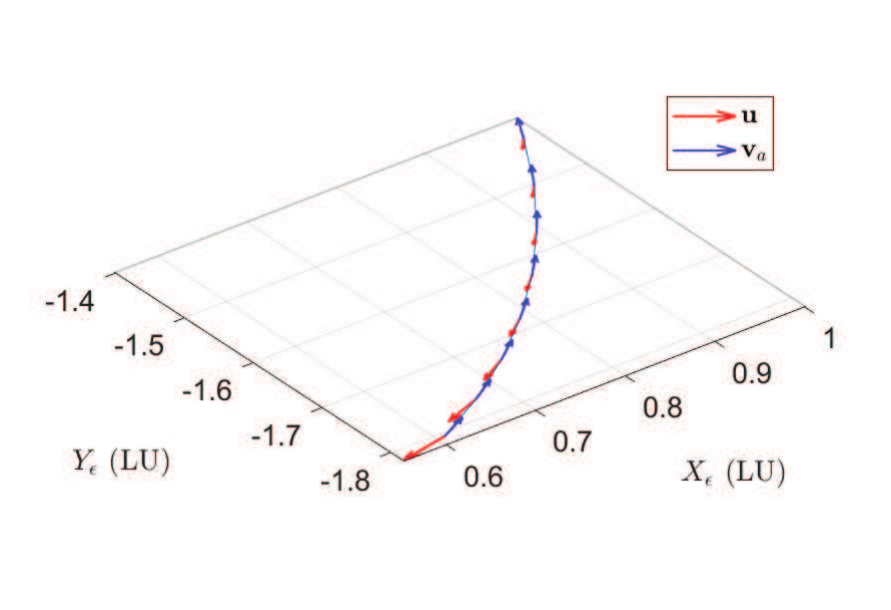}
	\captionof{figure}{Enlarge view of the operation region.}
	\label{fig:thesis_subopti_a12e06_sect_orbit_1034_lldec}
	\end{minipage}
\end{table}

\vspace*{-2em}

\section{Impulsive Control}\label{sec:4}
\begin{figure}[H]
	\subfloat[Case A]{
		\label{fig:thesis_impl_a11e02_orbit_casea}
		\includegraphics[width=0.5\textwidth]{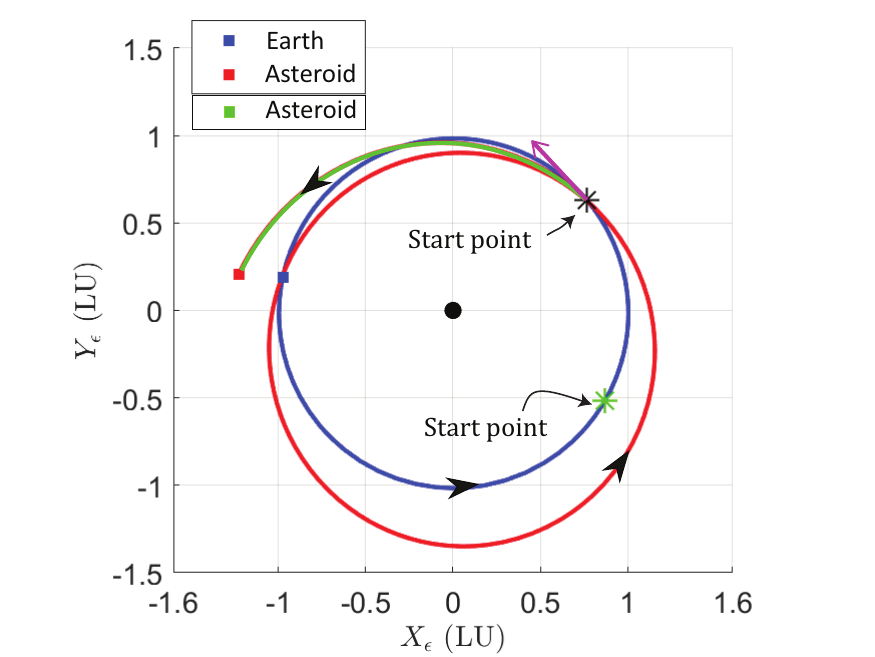} } 
	\subfloat[Case B]{
		\label{fig:thesis_impl_a11e02_orbit_caseb}
		\includegraphics[width=0.5\textwidth]{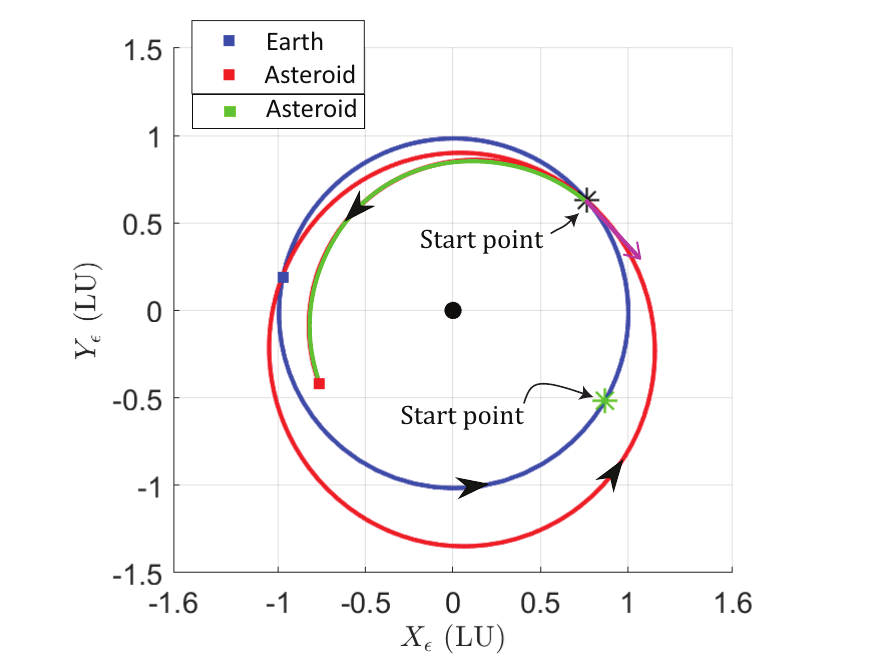} } 
	\caption{Showing orbits around the Sun --- The red colour shows the asteroid's orbit before the closest approach to Earth, and the green colour depicts the perturbed orbit after avoiding the collision with Earth.}
	\label{fig:thesis_impl_a11e02_orbit}
\end{figure}

As already mentioned in Section~\ref{sec:intro} Refs.~\cite{Ross2001} and~\cite{park_mazanek_2003} have calculated the optimal impulse change in velocity $\Delta \mathbf{v}$ and the authors have shown that there exist two optimal solutions to achieve the desired miss distance $\ell_m$ at any time. The two solutions have the same magnitude but have different phases, separated by $180$ degrees. However, these solutions lead to two different orbits for the ECO around the Sun after it misses the Earth. Due to the inclusion of the Earth's gravity on the ECO's trajectory, one solution modifies the ECO's orbit such that it gains momentum from the Earth during its close flyby. In contrast, the other solution leads to a decrease in momentum. The semi-major axis of the ECO's orbit increases in the former case whereas in the latter solution the semi major axis decreases. Figure \ref{fig:thesis_impl_a11e02_orbit} illustrates these results for an Apollo-type asteroid with orbital elements (similar to the asteroid Bennu)
\begin{table}[H]
 	\centering                            
 	\protect\caption{Orbital elements of the original orbit (Initial) }
 	\begin{tabular}{|l|l|l|l|}
 		\hline
 		a & e & i & Tp\\
 		$(au)$ &  & $(deg)$ &$(year)$\\
 		\hline 
 		1.1264 & 0.2037 & 0 & 1.1996\\
 		\hline 
 	\end{tabular}
 	\label{Tab:astelemin}
 \end{table}
The simulation is carried out for the impulse time equal to one period of asteroid ($t_{impulse} = 1\,Tp$). After the successful flyby, the orbit element become as listed in Table \ref{Tab:astelemout}. In this table, $T_p$ is the orbital period of asteroid and $\lambda$ is the impulse angle.
 \begin{table}[H]
	\centering                            
	\protect\caption{Orbital elements post flyby (Final).}
	\begin{tabular}{|l|l|l|l|l|l|l|l|l|l|}
		\hline
		\multicolumn{5}{|c|}{Case A} & \multicolumn{5}{|c|}{Case B} \\ 
		\hline
		$a$ & $e$ & $T_p$&$\Delta v$ & $\lambda$ & $a$ & $e$ &$T_p$& $\Delta v$ &$\lambda$ \\
		$(au)$ &  & $(year)$ &$ (cm/s) $ & $(deg)$ &$(au)$ &  & $(year)$ & $(cm/s) $ & $(deg)$\\
		\hline 
		1.3078 & 0.2729 & 1.5007& 58.90 &0 &0.9742 & 0.1812& 0.9649& 58.90& 180  \\
		\hline 
	\end{tabular}
	\label{Tab:astelemout}

\end{table}

\begin{figure}[H]
	\subfloat[Case A]{
		\label{fig:bennu_sep_excess}
		\includegraphics[width=0.5\textwidth]{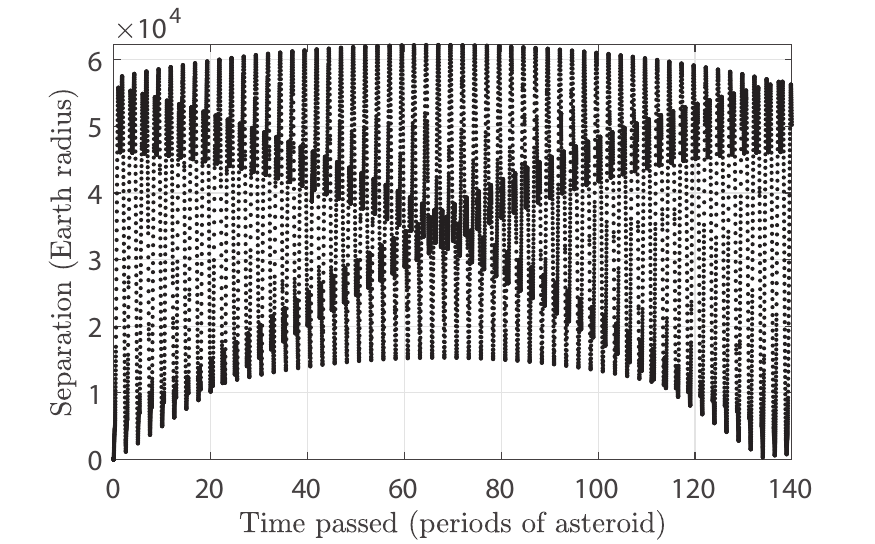} } 
	\subfloat[Case B]{
		\label{fig:bennu_sep_defi}
		\includegraphics[width=0.5\textwidth]{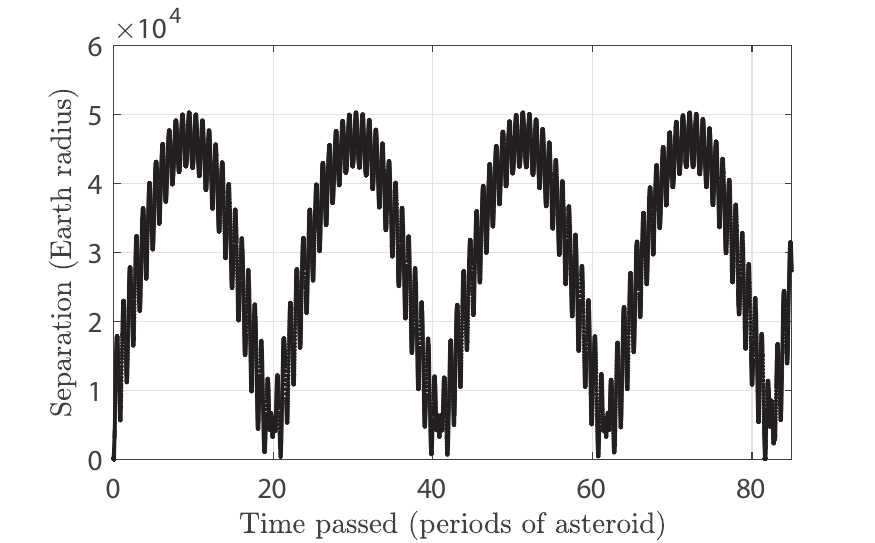} } 
	\caption{Variation of the distance between Earth and asteroid with time. The time instant $t=0$ corresponds to $\ell=\ell_{soi}$. }
	\label{fig:bennu_sep}
\end{figure}

In Figs.~\ref{fig:thesis_impl_a11e02_orbit} and \ref{fig:bennu_sep_excess}, we see that for Case A the impulse is applied in a direction making an acute angle ($\lambda=\,0$) with the ECO's velocity. As a result, the ECO crosses the point of intersection of two orbits after the Earth. In the process, the ECO gains momentum from Earth, and its total energy (in the Sun-ECO system) increases, increasing the heliocentric orbit's size. On the other hand, for the other impulse with $\lambda=180^{\circ}$, the ECO passes the point ahead of the Earth, reducing the ECO's momentum. Hence, its energy and the post flyby orbital size also decrease. The fact that the two solutions yield different trajectories for an asteroid around the Sun post flyby allows us to prefer one impulse over the other. The choice depends on the time interval for the ECO's next close pass by the Earth. In Case A, the ECO again passes after about 130 $T_p$, and in Case B, it takes around $20\,T_p$. For this example, Case A is more desirable than Case B. The same strategy can applied to results in continuous control and this will be done in a subsequent work.

\section{Conclusions}\label{sec:5}

This paper has designed optimal control strategies to deflect asteroids that cross Earth and has developed a general control strategy with variable power for operation. The variable power approach is found to provide a more energy-efficient control history in comparison with constant power operation
for the start point lying away from the perihelion. However, it is unable to utilize the additional time available before the pericenter pass. The constant power operation has yielded minimum operation time to achieve the necessary momentum change. However, the requirement of a constant power at all times made it less robust. 
It is advocated that the sub-optimal control provides viable solutions for the practical implementation as it requires less energy than the constant power method and eliminates the initial idle time. The time required by the sub-optimal solution lies between the constant power operation and variable power operation. The control history is easier to implement from a practical viewpoint. Also, the best time to start operation for the variable power approach is the earliest possible perihelion pass. A comparison is done for the impulsive transfer methods. Among the two solutions provided by the existing work, the finding of this work has shown that one solution is more favorable as it increases the time interval until a future crossing of the same asteroid with the Earth. Comparing the constant power, variable power, and impulsive thrust, the study identified situations when one control strategy could be better than the others.   

Finally, the paper presents the effect of the Moon on the asteroid-Earth collision events (see Appendix). Due to the short time for the transit, Moon was found to exhibit extremely less influence on the orbit of the asteroid. However, depending on the timing and phasing of the approach, Moon was found to either pull the asteroid more towards Earth or push away from Earth. This can be useful in the controlled asteroid capture process for mining applications.

\section*{Appendix: Effect of Moon's Gravity}\label{sec:appn}
\begin{figure}[H]
	\centering
	\includegraphics[width= 0.6\textwidth] {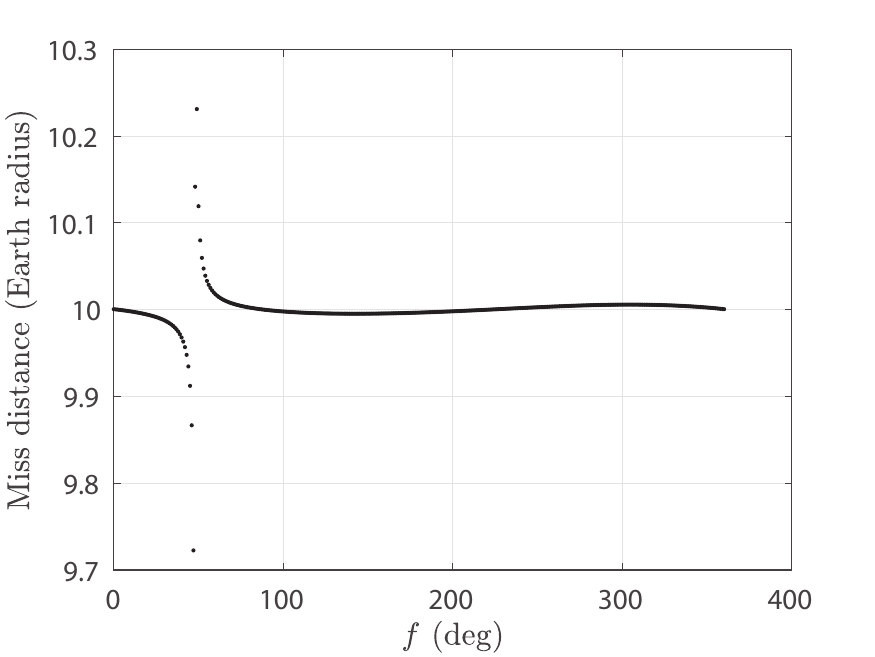}	  
	\caption{The graph shows variation of the miss distance with true anomaly of Moon.}
	\label{fig:thesis_moon_rm_vs_f}
\end{figure}

\begin{figure}[H]
	\hspace*{-3em}
	\subfloat[For $f=47^\circ$]{
		\label{fig:mooneffect_m}
		\includegraphics[width=0.5\textwidth]{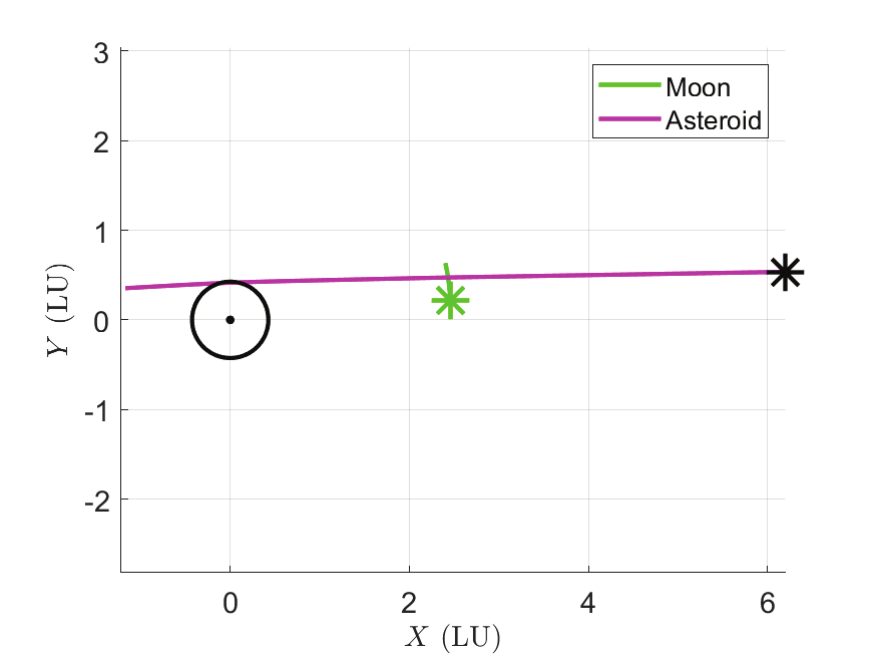} } 
	\hspace*{-3em}
	\subfloat[For $f=48^\circ$]{
		\label{fig:mooneffect_p}
		\includegraphics[width=0.5\textwidth]{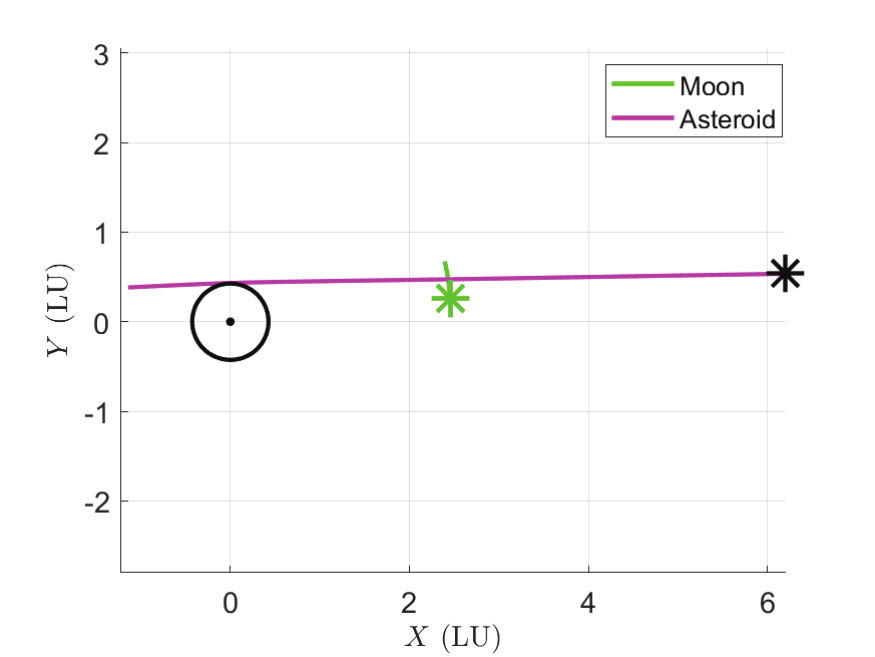} } 
	\caption{The graph shows the orbits of asteroid and Moon around Earth. The asteroid starts from the edge of Earth's SOI (see far right in the graph indicated by a black asterisk). The black circle has a radius of 10 Earth radius. (A) shows orbit orientation for maximum reduction in $\ell_m$, and (B) shows orbit orientation for maximum gain in $\ell_m$.}
	\label{fig:mooneffect_Plus and minus}
\end{figure}

In this appendix, we investigate the effect of Moon’s gravity on the miss distance for an asteroid inside the SOI of Earth. The objective is to find out the change in \textit{minimum separation} between Earth and asteroid in the presence of Moon. Consider a scenario when the asteroid is at the boundary of SOI and the Moon is at some true anomaly $f$ in its orbit around Earth. This configuration will serve as the initial condition for the analysis. For simplicity, the inclination\footnote{The actual inclination of the Moon’s orbit relative to the ecliptic plane is 5 degrees} of the Moon orbit is considered to be zero.
For the simulation, the original orbital elements of the asteroid’s orbit as specified in Table \ref{tab:ast_moon} are considered. The nominal miss distance of the asteroid (without including the Moon’s gravity) is 10 Earth radius. We will find the deviation of the  miss distance from the nominal value for the various initial Earth-Moon positional configurations. Numerical scaling of variables --- the length unit $(LU)$, $1\,LU = 1\, au$, time unit $(TU)$, $1\,TU= 1\ day$, and the speed unit becomes $(SU)$, $1\,SU = 1\,au/1\,day$.
The simulation shows that the relative error for the lowest miss distance is only 2.77 \% for normal miss distance $\ell_m = 10$ Earth radii (nominal miss distance refers to the case when Moon’s gravity is not present), refer Figs.~\ref{fig:thesis_moon_rm_vs_f}--\ref{fig:mooneffect_Plus and minus}. Thus, we can conclude that neglecting Moon gravity has not affected the above results significantly. However, one needs to be careful if the desired miss distance is close to 1 Earth radius as it can be for larger asteroids to minimise fuel consumption. Also, asteroids can tidally break due to moon's gravity. Hence, it is essential to avoid both Earth and moon to avoid creating tidally disrupted asteroid fragments that can pose threats in future. At outcome of such an event is also a subject for a future study.
\begin{table}[h]
	\centering
	\caption{Specification of the asteroid's orbit around the Sun and Moon's orbit around the Earth}
	\begin{subtable}[h]{0.3\textwidth}
		\begin{tabular}{| m{4em} | m{2em} | m{4em} | }
			\hline
			a ($au$) & e & i ($deg$)\\
			\hline 
			1.2 & 0.6 & 0\\
			\hline 
		\end{tabular}
		\caption{Asteroid }
		\label{tab:ast_moon}
	\end{subtable}
\hspace*{3cm}
	\begin{subtable}[h]{0.3\textwidth}
		\begin{tabular}{| m{4em} | m{3em} | m{3em} |}
			\hline
			a ($km$)& e & i ($deg$) \\
			\hline 
			$0.3844e6$  & 0.0549 & 0\\
			\hline
		\end{tabular}
		\caption{Moon }
		\label{Tab:moonandast}
	\end{subtable}
	\label{tab:Moon}
\end{table}
\vspace{-2em}

\bibliography{main}

\begin{thebibliography}{20}
\newcommand{\enquote}[1]{``#1''}
\providecommand{\natexlab}[1]{#1}
\providecommand{\url}[1]{\texttt{#1}}
\providecommand{\urlprefix}{URL }
\expandafter\ifx\csname urlstyle\endcsname\relax
  \providecommand{\doi}[1]{\discretionary{}{}{}https://doi.org/#1}\else
  \providecommand{\doi}[1]{\discretionary{}{}{}\urlstyle{rm}\url{https://doi.org/#1}}\fi

\bibitem[{Solem(1993)}]{doi:10.2514/3.11531}
Solem, J.~C., \enquote{Interception of comets and asteroids on collision course
  with Earth,} \emph{Journal of Spacecraft and Rockets}, Vol.~30, No.~2, 1993,
  pp. 222--228.
\newblock \doi{10.2514/3.11531},
  \urlprefix\url{https://doi.org/10.2514/3.11531}.

\bibitem[{Ahrens and Harris(1992)}]{ahrens_harris_1992}
Ahrens, T.~J., and Harris, A.~W., \enquote{Deflection and fragmentation of
  near-earth asteroids,} \emph{Nature}, Vol. 360, No. 6403, 1992, p. 429–433.
\newblock \doi{10.1038/360429a0}.

\bibitem[{Park and Mazanek(2003)}]{park_mazanek_2003}
Park, S.-Y., and Mazanek, D.~D., \enquote{Mission Functionality for Deflecting
  Earth-Crossing Asteroids/Comets,} \emph{Journal of Guidance, Control, and
  Dynamics}, Vol.~26, No.~5, 2003, pp. 734--742.
\newblock \doi{10.2514/2.5128}, \urlprefix\url{https://doi.org/10.2514/2.5128}.

\bibitem[{McInnes(1999)}]{McInnes1999}
McInnes, C.~R., \emph{Solar sail orbital dynamics}, Springer London, London,
  1999, pp. 112--170.
\newblock \doi{10.1007/978-1-4471-3992-8_4},
  \urlprefix\url{https://doi.org/10.1007/978-1-4471-3992-8_4}.

\bibitem[{Maddock et~al.(2007)Maddock, Sánchez, Vasile, and
  Radice}]{articleMaddock}
Maddock, C., Sánchez, J.-P., Vasile, M., and Radice, G., \enquote{Comparison
  of Single and Multi-Spacecraft Configurations for NEA Deflection by Solar
  Sublimation,} \emph{AIP Conference Proceedings}, Vol. 886, 2007, pp.
  303--316.
\newblock \doi{10.1063/1.2710064}.

\bibitem[{Park and Ross(1999)}]{ParkandRoss1999}
Park, S., and Ross, I., \enquote{Two-body optimization for deflecting
  earth-crossing asteroids,} \emph{Journal of Guidance, Control, and Dynamics},
  Vol.~22, No.~3, 1999, pp. 415--420.
\newblock \doi{10.2514/2.4413}.

\bibitem[{Song et~al.(2010)Song, Park, and Choi}]{articleSong2010}
Song, Y.-J., Park, S.-Y., and Choi, K.-H., \enquote{Mission feasibility
  analysis on deflecting Earth-crossing objects using a power limited laser
  ablating spacecraft,} \emph{Advances in Space Research}, Vol.~45, No.~1,
  2010, pp. 123--143.
\newblock \doi{https://doi.org/10.1016/j.asr.2009.08.023},
  \urlprefix\url{https://www.sciencedirect.com/science/article/pii/S0273117709005961}.

\bibitem[{Ross et~al.(2001)Ross, Park, and Porter}]{Ross2001}
Ross, I., Park, S., and Porter, S., \enquote{Gravitational effects of earth in
  optimizing delta-V for deflecting earth-crossing asteroids,} \emph{Journal of
  Spacecraft and Rockets}, Vol.~38, No.~5, 2001, pp. 759--764.
\newblock \doi{10.2514/2.3743}.

\bibitem[{Conway(2001)}]{doi:10.2514/2.4814}
Conway, B.~A., \enquote{Near-Optimal Deflection of Earth-Approaching
  Asteroids,} \emph{Journal of Guidance, Control, and Dynamics}, Vol.~24,
  No.~5, 2001, pp. 1035--1037.
\newblock \doi{10.2514/2.4814}, \urlprefix\url{https://doi.org/10.2514/2.4814}.

\bibitem[{{Bekey}(1997)}]{1997ESASP.393..699B}
{Bekey}, I., \enquote{{Project Orion: Orbital Debris Removal Using Ground-Based
  Sensors and Lasers},} \emph{Second European Conference on Space Debris}, ESA
  Special Publication, Vol. 393, edited by B.~{Kaldeich-Schuermann}, 1997, p.
  699.
\newblock \doi{10.2514/6.1997-630}.

\bibitem[{Phipps(2003)}]{articlephipps2003}
Phipps, C., \enquote{Overview of laser applications: The state of the art and
  the future trend,} \emph{Proceedings of SPIE - The International Society for
  Optical Engineering}, Vol. 4831, 2003.
\newblock \doi{10.1117/12.486517}.

\bibitem[{Song and Park(2015)}]{articleSong2015}
Song, Y.-J., and Park, S.-Y., \enquote{Estimation of necessary laser power to
  deflect near-Earth asteroid using conceptual variable-laser-power ablation,}
  \emph{Aerospace Science and Technology}, Vol.~43, 2015, pp. 165--175.
\newblock \doi{10.1016/j.ast.2015.02.022}.

\bibitem[{Chang-Diaz et~al.(2000)Chang-Diaz, Squire, Bengtson, Breizman,
  Carter, and Baity}]{ChangDiaz2000}
Chang-Diaz, F., Squire, J., Bengtson, R., Breizman, B., Carter, M., and Baity,
  F., \emph{The physics and engineering of the VASIMR engine}, 2000.
\newblock \doi{10.2514/6.2000-3756},
  \urlprefix\url{https://arc.aiaa.org/doi/abs/10.2514/6.2000-3756}.

\bibitem[{Park and Choi(2005)}]{ParkYoung2005}
Park, S.-Y., and Choi, K.-H., \enquote{Optimal Low-Thrust Intercept/Rendezvous
  Trajectories to Earth-Crossing Objects,} \emph{Journal of Guidance, Control,
  and Dynamics}, Vol.~28, No.~5, 2005, pp. 1049--1055.
\newblock \doi{10.2514/1.10895},
  \urlprefix\url{https://doi.org/10.2514/1.10895}.

\bibitem[{Gambi et~al.(2022)Gambi, {García del Pino}, Mosser, and
  Weinmüller}]{GAMBI2021}
Gambi, J., {García del Pino}, M., Mosser, J., and Weinmüller, E.,
  \enquote{Trailing formations of lightweight spacecrafts to deflect NEAs by
  means of laser ablation,} \emph{Acta Astronautica}, Vol. 190, 2022, pp.
  241--250.
\newblock \doi{https://doi.org/10.1016/j.actaastro.2021.10.006},
  \urlprefix\url{https://www.sciencedirect.com/science/article/pii/S0094576521005476}.

\bibitem[{{Wikipedia contributors}(2022)}]{wikiLaser}
{Wikipedia contributors}, \enquote{Laser --- {Wikipedia}{,} The Free
  Encyclopedia,} , 2022.
\newblock
  \urlprefix\url{https://en.wikipedia.org/w/index.php?title=Laser&oldid=1071093900},
  [Online; accessed 2-March-2022].

\bibitem[{Phipps(2007)}]{phipps2007laser}
Phipps, C.~R., \emph{Laser ablation and its applications}, Springer, 2007, pp.
  375--549.
\newblock \doi{10.1007/978-0-387-30453-3}.

\bibitem[{Battin et~al.(1987)Battin, of~Aeronautics, and
  Astronautics}]{battin1987introduction}
Battin, R., of~Aeronautics, A.~I., and Astronautics, \emph{An Introduction to
  the Mathematics and Methods of Astrodynamics}, AIAA Textbook Series, American
  Institute of Aeronautics and Astronautics, 1987, Chap.~9.
\newblock \urlprefix\url{https://books.google.co.in/books?id=xZBTAAAAMAAJ}.

\bibitem[{bet(2010)}]{betts2010practicalch4}
\emph{Practical Methods for Optimal Control and Estimation Using Nonlinear
  Programming, Second Edition}, Society for Industrial and Applied Mathematics
  (SIAM), 2010, Chap.~4, pp. 123--217.
\newblock \doi{10.1137/1.9780898718577.ch4},
  \urlprefix\url{https://epubs.siam.org/doi/abs/10.1137/1.9780898718577.ch4}.

\bibitem[{Gibbings et~al.(2013)Gibbings, Vasile, Watson, Hopkins, and
  Burns}]{GIBBINGS2013}
Gibbings, A., Vasile, M., Watson, I., Hopkins, J.-M., and Burns, D.,
  \enquote{Experimental analysis of laser ablated plumes for asteroid
  deflection and exploitation,} \emph{Acta Astronautica}, Vol.~90, No.~1, 2013,
  pp. 85--97.
\newblock \doi{https://doi.org/10.1016/j.actaastro.2012.07.008},
  \urlprefix\url{https://www.sciencedirect.com/science/article/pii/S0094576512002718},
  nEO Planetary Defense: From Threat to Action - Selected Papers from the 2011
  IAA Planetary Defense Conference.

\end{thebibliography}

\end{document}